\documentclass[reqno]{amsart}

\input xy
\xyoption{all}
\usepackage{accents}
\usepackage{epsfig}
\usepackage{color}
\usepackage{amsthm}
\usepackage{amssymb}
\usepackage{amsmath}
\usepackage{amscd}
\usepackage{amsopn}
\usepackage{graphicx}
\usepackage{centernot}
\usepackage{xspace}

\usepackage{url}
\usepackage{enumitem, hyperref}\hypersetup{colorlinks}

\usepackage{color} 

\definecolor{darkred}{rgb}{1,0,0} 
\definecolor{darkgreen}{rgb}{0,0.6,0}
\definecolor{darkblue}{rgb}{0,0,.8}

\hypersetup{colorlinks,
linkcolor=darkblue,
filecolor=darkgreen,
urlcolor=darkred,
citecolor=darkgreen}

\makeatletter
\def\reflb#1#2{\begingroup
    #2%
    \def\@currentlabel{#2}%
    \phantomsection\label{#1}\endgroup
}
\makeatother

\numberwithin{equation}{section}
\newtheorem {Theorem}{Theorem}
\numberwithin{Theorem}{section}
\newtheorem{Sublemma}{Lemma}[Theorem]

\newtheorem {Lemma}[Theorem]    {Lemma}

\newtheorem {Proposition}[Theorem]{Proposition}

\theoremstyle{definition}
\newtheorem{Definition}[Theorem]{Definition}
\theoremstyle{remark}
\newtheorem{Remark}[Theorem]{Remark}
\newtheorem{Example}[Theorem]{Example}

\def    \12    {{\frac{1}{2}}}

\newcommand{\overbar}[1]{\mkern 1.5mu\overline{\mkern-1.5mu#1\mkern-1.5mu}\mkern 1.5mu}

\usepackage{mathtools}

\usepackage{tikz}
\usepackage{tikz-cd}

\usepackage{xspace}
\usepackage[foot]{amsaddr}

\begin{document}


\setlength{\smallskipamount}{6pt}
\setlength{\medskipamount}{10pt}
\setlength{\bigskipamount}{16pt}





\title[Extending the Mean Index]{On an Extension of the Mean Index to a large subset of Linear Canonical Relations}
\address{MG: Department of Mathematics, UC Santa Cruz, Santa Cruz, CA 95064, USA} 
\author[Yusuf G\"{o}ren]{Yusuf G\"{o}ren$^\dag$}
\author[Matt Grace]{Matthew Grace}
\thanks{$^\dag$ The author's contributions to the paper occurred before current employment at Amazon}
\thanks{The work is partially supported by NSF grant DMS-1308501}

\email{yusuf.goren@gmail.com, migrace@ucsc.edu}



\date{\today} 

\bigskip

\begin{abstract}
	
	In this paper, viewing the symplectic linear group as a subset of the Lagrangian Grassmannian we extend the mean index to the complement of a codimension-two subset of the Grassmannian. This extension retains many of the desirable properties of the mean index, the most significant of which are continuity and a homogeneity condition adapted to the set-theoretic composition of canonical relations.
\end{abstract}
\maketitle

\tableofcontents

\section{Introduction}
\label{sec:intro+results}
In this paper we prove the existence of a continuous extension $\hat{\Delta}$ of the mean index $\Delta$ via the identification of the symplectic group with an open dense subset of the Lagrangian Grassmannian (definition~\ref{def:graph}). The mean index $\Delta$ may be constructed as a real valued map defined over arbitrary paths in $Sp(2n)$ although this is not the most useful interpretation. Most applications of the mean index tend to restrict $\Delta$ to paths originating at the identity as the collection of all such paths in $Sp(2n)$ may be identified up to homotopy with $\widetilde{Sp}(2n)$, the universal cover of the symplectic group. This paper will use the former notion (remark~\ref{rmk:fund} provides motivation for this choice) to construct $\hat{\Delta}$ as a real valued map defined over all paths in the Lagrangian Grassmannian which never intersect a certain set $H$ of codimension two in the Grassmannian. The bulk of the proof is in showing the existence of a continuous, circle valued map $\hat{\rho}$ defined on the complement of $H$ in the Grassmannian which continuously extends the square (see remark~\ref{rmk:squared}) of the circle map $\rho\colon Sp(2n) \to S^1$, perhaps most notable for coinciding with the complex determinant when restricted to the unitary group $U(n) \cong Sp(2n) \cap O(2n)$. Indeed since $\rho$ gives rise to $\Delta$ via a formal construction (see definition~\ref{def:mean}) the existence and continuity of $\hat{\Delta}$ is an immediate consequence of the continuity of $\hat{\rho}$ when constructed in the same manner (an early example of this lifting procedure may be found in \cite{Mi}). The extension is shown to preserve some of the desirable qualities of the mean index, of which fixed endpoint homotopy invariance and a homogeneity condition adapted to the set-theoretic composition of canonical relations are most relevant.


The mean index may alternatively be defined (see \cite{Lon08} for a thorough treatment) for any path $\gamma \subset Sp(2n)$ originating at the identity via a generalized Conley-Zehnder index (see remark~\ref{rmk:fund} for the definition of the composition of paths originating at the identity) as follows, 
\begin{equation}
 \Delta(\gamma) := \lim_{k \to \infty} \frac{\mu_{cz}(\gamma^k)}{k}.\nonumber
\end{equation} A critical fact relating the two indices is the following bound which holds for paths in $Sp(2n)$ originating at the identity, 
\begin{equation}\label{eqn:bound}
|\Delta(\gamma) - \mu_{cz}(\gamma)| \leq n.
\end{equation} An early example of this bound is the strict inequality found in \cite{SaZ} holding for all `admissible' paths $\gamma$. This was later extended to the above in which degenerate $\gamma$ are the only paths for which equality may occur. 

In addition to satisfying ~\eqref{eqn:bound}, the mean index possesses other algebraic properties, of which the fact that $\Delta$ is the unique, real valued homogeneous quasimorphism on $\widetilde{Sp}(2n)$ \cite{BG} is crucial. These properties have yielded results \cite{GG15,GG09} in counting Reeb orbits as well as identifying the multiplicity of the periodic orbits of a Hamiltonian system within the framework of various Floer homologies. One of the seminal results in this regard is found in the aforementioned \cite{SaZ} in which the existence of infinitely many periodic orbits of a given Hamiltonian system is shown provided some non-degeneracy conditions are held by the Hamiltonian at the orbits, i.e. the Conley conjecture. An example of advances in this direction may be found in \cite{GG10} in which a generalization of the Conley conjecture is proven via local Floer homology wherein any Hamiltonian diffeomorphism of a closed, symplectically aspherical manifold with isolated fixed points possesses simple periodic points of arbitrarily large period. A proof of the classical Conley Conjecture for closed, symplectically aspherical manifolds utilizing the mean index and the resulting filtration of the local Floer algebra of an isolated periodic point may be found in \cite{Ci}. The index theory common to all of these references may be found in \cite{GG15} wherein the framework of Lusternik-Schnirelmann theory is applied to various symplectic and contact homology theories.

The price paid in extending the mean index from the symplectic group to the complement of $H$ in the Lagrangian Grassmannian is perhaps most plainly demonstrated by the lack of group structure over the linear canonical relations of a fixed symplectic vector space when equipped with the set-theoretic composition operation (definition~\ref{def:comp}). In particular this implies the quasimorphism property $\Delta$ possesses has no immediate analogue for $\hat{\Delta}$ in the context of linear canonical relations and set-theoretic composition. The composition map for linear canonical relations is even discontinuous with respect to the usual topology on the Lagrangian Grassmannian, a fact which is compounded by the presence of severe complications involved in transitioning categorically to smooth canonical relations \cite{Wei}. Regardless, the construction of $\hat{\Delta}$ for linear canonical relations may be a step towards developing machinery allowing for the further study of Reeb orbits and the evolution and periodicity of Hamiltonian systems within the framework of canonical relations. 


\medskip
\noindent\textbf{Acknowledgments} The authors would like to thank Viktor Ginzburg for his helpful advice and guidance throughout the writing of this paper as well as Alan Weinstein for our correspondences.

\section{Definitions and Conventions}
\label{sec:prelim}
Let $(V^{2n},\omega)$ be a symplectic vector space and let $V \times \overbar{V}$ denote $(V\times V,\tilde{\omega} = \pi_1^\ast\omega - \pi_2^\ast\omega)$. A Lagrangian subspace $L \leq V \times \overbar{V}$, also referred to as a Lagrangian (canonical) relation is said to have source $V$ and target $\overbar{V}$ although many authors define the target and source backwards relative to this definition. We introduce the following notation used in \cite{To} for arbitrary linear relations (with the term $halo(L)$ introduced in \cite{LorW} and more recently denoted $indet(L)$ as in \cite{LiW}).
\begin{Definition} 
	
	Given any linear canonical relation $L \in \Lambda_{2n} := LagGr(V\times \overbar{V},\tilde{\omega})$ denote the following distinguished subspaces of $V$.
\begin{itemize}
	\item $ker(L) := \left\{v \in V \, | \, (v,0) \in L\right\}$
	\item $halo(L) := \left\{v \in V \, | \, (0,v) \in L\right\}$
	\item $dom(L) := \left\{v \in V \, | \, \exists w \in V \, , \, (v,w) \in L\right\} = ker(L)^\omega$
	\item $ran(L) := \left\{v \in V \, | \, \exists w \in V \, , \, (w,v) \in L\right\} = halo(L)^\omega$
\end{itemize}
\end{Definition}
It may be shown that both $dom(L)$ and $ran(L)$ are co-isotropic and share the same dimension since the target is equal to the source: $(V,\omega)$.

\begin{Definition}\label{def:graph}
	Define the following smooth map sending each $A \in Sp(V)$ to its graph, a Lagrangian subspace of $V \times \overbar{V}$.
\begin{eqnarray}
Sp(V) &\underset{Gr}{\hookrightarrow}& \Lambda_{2n}\nonumber\\
A &\mapsto& \left\{(v,Av) \in V\times\overbar{V} \, | \, v \in V \right\}\nonumber 
\end{eqnarray} 
\end{Definition}
The above map has an open and dense image in $\Lambda_{2n}$ \cite{He} so that $\Lambda_{2n}$ is a compactification of $Sp(2n)$.

The continuous circle map $\rho\colon Sp(V) \to S^1$, perhaps most notable for extending the circle valued complex determinant defined on the subgroup $U(n) \cong Sp(2n,\mathbb{R}) \cap O(2n,\mathbb{R})$ lies at the heart of the mean index.
\begin{Definition}\label{def:mean}
given any path $\gamma\colon I \to Sp(V)$ there exists a unique, continuous $\theta\colon I \to \mathbb{R}$ such that $(\rho\circ \gamma)(t) = e^{i\theta(t)}$ and $\theta(0) \in [-\pi,\pi)$. Then the mean index for the path $\gamma$ is defined as \begin{eqnarray}
\Delta(\gamma) := \frac{\theta(1)-\theta(0)}{2\pi}\nonumber.
\end{eqnarray} 
\end{Definition}

\begin{Remark}\label{rmk:squared}
	The usual definition of the mean index will need to be altered before attempting to extend $\rho$ as shown below in example~\ref{example:squared}, we will solve this by continuously extending $\rho^2$ instead. Despite this difference the extended mean index $\hat{\Delta}$ may be shown to differ from $\Delta$ by a factor of two.
\end{Remark}

\begin{Remark}\label{rmk:fund}
	A critical property of the mean index as used in the papers referenced in the introduction relies on the fact that $\pi_1(Sp(2n)) \cong \mathbb{Z} \cong \pi_1(S^1)$ and that $\rho$ is a circle map (that is, the induced map for $\rho$ yields an isomorphism of fundamental groups). This allows homotopy classes of paths originating at the identity to be identified uniquely with $g \in \widetilde{Sp}(2n)$ meaning that the mean index then becomes the unique continuous, homogeneous quasimorphism $\Delta\colon \widetilde{Sp}(2n) \to \mathbb{R}$ where the composite of two paths is defined using the group structure of $\widetilde{Sp}(2n)$ and the above identification. In example~\ref{example:fundamental} we consider for $n=1$ how to partially address the obstruction $\pi_1(\Lambda_{2n} \setminus H)$ presents in reconstructing this process with canonical relations.
\end{Remark}

\begin{Definition}\label{def:comp}
	The set theoretic composition for linear (canonical) correspondences is defined below given $L,L^\prime \in \Lambda_{2n}$,
	\begin{eqnarray}
	L \circ L^\prime = \left\{(x,z) \, | \, \,\, \exists y \in ran(L^\prime)\cap dom(L) \text{ such that } (x,y) \in L^\prime, \,\, (y,z) \in L \right\}.\nonumber
	\end{eqnarray}
\end{Definition}

As first formulated in \cite{Ho} following \cite{Ma}, it is necessary to impose the transversality condition $dom(L) \pitchfork ran(L^\prime)$ since the above set theoretic operation, while well defined, fails to be continuous without it (see example~\ref{example:cont}). The category of canonical relations is often denoted \textbf{SLREL}, and (to the best of the author's knowledge) was first formally constructed in \cite{BT}. A diverse collection of techniques have been utilized to avoid first, the obstacle of discontinuity in the linear case and second, the myriad of additional problems that arise in the smooth case; one early solution may be found in \cite{GuiS} in which the authors augment Lagrangian relations with half densities within the context of microlocal analysis of certain integral transforms. This technique and others are mentioned in \cite{Wei} which approaches the aforementioned issues from a categorical perspective. We postpone to section~\ref{subsec:smooth} both the many issues in adapting the extended mean index to smooth objects as well as the introduction of the most promising `extension' of canonical relations for our purposes, first developed in \cite{WehW} and later refined in \cite{LiW}.
\section{Main Theorems and Proof Breakdown}\label{subsec:sketch}

\begin{Theorem}\label{thm:mainthm}
	Define $H := \left\{L \in \Lambda_{2n}\, | \, dim(dom(L) \cap halo(L)) \ge 1\right\}$,
	
	\begin{enumerate}
		
		\item Theorem~\ref{thm:codim}: The codimension of $H$ is two.
		\item Theorem~\ref{subthm:mainthm:2}: There exists a unique, real valued continuous function $\hat{\Delta}$ constant on fixed endpoint homotopy classes of paths in $\Lambda_{2n}\setminus H$ such that for any path $\gamma\colon I\to Sp(2n)$ we have that $\hat{\Delta}(Gr(\gamma)) = 2\Delta(\gamma)$. 
		\item  Theorem~\ref{subthm:mainthm:3}: The extended mean index $\hat{\Delta}$ is homogeneous; given any $\gamma\colon I \to \Lambda_{2n}\setminus H$ with $(\gamma \cdot \gamma)(t) := \gamma(t) \circ \gamma(t)$ (where $\circ$ is defined as in definition~\ref{def:comp} above) then $\hat{\Delta}(\gamma^k) = k\cdot \hat{\Delta}(\gamma)$ for all $k \ge 0$.
	\end{enumerate}
\end{Theorem}

The most intricate claim above is part 2 and the key to proving it is showing that $\rho^2$ may be extended continuously to $\hat{\rho}\colon\Lambda_{2n}\setminus H \to S^1$ (Theorem~\ref{thm:rhothm}), after which $\hat{\Delta}$ may be constructed. Indeed, provided a continuous extension $\hat{\rho}$ exists we have for any path $\gamma\colon I \to \Lambda_{2n}\setminus H$ a unique continuous $\hat{\theta}\colon I \to \mathbb{R}$ such that $(\hat{\rho}\circ \gamma)(t) = e^{i\hat{\theta}(t)}$ and $\hat{\theta}(0) \in [-\pi,\pi)$. Then the extended mean index for the path $\gamma$ may be defined as $\hat{\Delta}(\gamma) := \frac{\hat{\theta}(1)-\hat{\theta}(0)}{2\pi}$ so that part 2 of theorem~\ref{thm:mainthm} is true if theorem~\ref{thm:rhothm} is. Note that the continuity and homotopy invariance of $\hat{\Delta}$ is a formal consequence of the continuity of $\hat{\rho}$. $\qed$

The exceptional set $H$ manifests in the two dimensional case as a circle bridging the two connected components of the parabolic transformations at a projective `line at infinity' outside the image $Gr(Sp(2)) \subset \Lambda_2$ (see figure~\ref{fig:graph}). The fact that $\hat{\rho}$ is not a circle function is in some sense a consequence of removing $H$ from $\Lambda_{2n}$; if $A_i \in Sp(2)$ with $Gr(A_i) \to L \in H \cong \left\{\langle v\rangle \times \left\{0\right\} \oplus \left\{0\right\} \times \langle v \rangle  | \,\, v \in V\right\}$ then the value of $\lim_{i\to\infty}\rho(A_i)$ depends in particular on whether one is eventually approaching via elliptic or hyperbolic sequences. If $Tr(A_i) > 2$ for all $i \in \mathbb{N}$ and $Gr(A_i) \to L$ then $\rho(A_i) = \pm 1$ for all $i$ so that $lim_{i\to\infty}\rho(A_i) = \pm 1$. On the other hand when approaching $H$ from the elliptic transformations significant discontinuities arise even in low dimensions; in \cite{G15} a family of sequences of symplectic maps $\left\{A_i^\theta\right\} \in Sp(2)$ are constructed for each $\theta \in S^1\setminus \left\{\pm 1\right\}$ such that for a fixed $L \in H$ each  $\lim_{i \to \infty}Gr(A_i^\theta) = L$ yet $\rho(A_i^\theta) = \theta$ for all $i \in \mathbb{N}$.

Now we state the theorem at the core of the proof for part 2 of theorem~\ref{thm:mainthm}.

\begin{Theorem}\label{thm:rhothm}
	There exists a unique continuous map $\hat{\rho}\colon \Lambda_{2n}\setminus H \to S^1$ such that for any $L \in \Lambda_{2n}\setminus H$ satisfying $L = Gr(\phi)$ for some $\phi \in Sp(V)$ then $\hat{\rho}(Gr(\phi)) = \rho^2(\phi)$.
\end{Theorem}

To extend $\rho^2$ will require proving the following three claims:

\begin{enumerate}[wide, labelwidth=!, labelindent=0pt]
 \item Theorem~\ref{subthm:decomp}: For a given $L \in \Lambda_{2n}\setminus H$ there exists a unique symplectic decomposition of $V = V_s \oplus V_g$ and $\phi \in Sp(V_g)$ such that 
\begin{eqnarray}\label{eqn:decomp}
L = (ker(L) \times \left\{0\right\}) \oplus (\left\{0\right\} \times halo(L)) \oplus Gr(\phi)\leq (V_s\times \overbar{V}_s) \oplus (V_g \times \overbar{V}_g)
\end{eqnarray}
where $ker(L),halo(L) \in LagGr(V_s)$ are transverse.		

\item Theorem~\ref{subthm:hyperbolic}: Given any sequence $\left\{A_i\right\}_{i=1}^\infty \subset Sp(V)$ where each $A_i$ has distinct eigenvalues such that 
\begin{eqnarray}
Gr(A_i) \to L_1 \times \left\{0\right\} \oplus \left\{0\right\} \times L_2 \in \Lambda_{2n}\setminus H\nonumber
\end{eqnarray} 
(i.e. $L_1,L_2 \in \Lambda_n$ are transversal), then the $A_i$ will eventually have no eigenvalues $\lambda \in S^1\setminus \left\{\pm 1\right\}$. In particular this shows that $\rho^2(A_i) = 1$ for sufficiently large $i$.
	
\item Theorem~\ref{subthm:rhothm:3}: Any sequence $\left\{A_i\right\}_{i=1}^\infty  \subset Sp(V)$ each with distinct eigenvalues which has 
	\begin{eqnarray}
	Gr(A_i) \underset{i\to\infty}{\to} L \in \Lambda_{2n}\setminus H\nonumber
	\end{eqnarray} 
	such that the graph part $\phi$ of $L$ (see ~\eqref{eqn:decomp}) has semisimple eigenvalues induces a sequence of unique $A_i$ invariant symplectic decompositions $V = E^i_s \oplus E^i_g$ so that $A_i = \alpha_i\oplus \beta_i \in Sp(E^i_s) \times Sp(E^i_g)$ which satisfies $Gr(\alpha_i) \to ker(L) \times \left\{0\right\} \oplus \left\{0\right\} \oplus halo(L)$. 
	
	Additionally there exists an $N \in \mathbb{N}$ for which there is a sequence of symplectic isomorphisms
	\begin{eqnarray}
	\left\{I_i\colon (E^i_g,\left.\omega\right|_{E^i_g \times E^i_g})  \overset{\cong}{\to} (V_g,\left.\omega\right|_{V_g \times V_g}) \right\}_{i=N}^\infty\nonumber
	\end{eqnarray} uniquely determined by $L$ such that each $\beta_i\colon E^i_g \to E^i_g$ is conjugate via $I_i$ to some $\phi_i \in Sp(V_g)$ for all $i \ge N$ with $\phi_i \to \phi$. Most importantly we show this preserves the data used in computing $\rho$, namely the eigenvalues and the conjugacy classes of the $A_i$ restricted to elliptic eigenspaces.
\end{enumerate}
	\begin{Remark}
	Refer to the original paper \cite{SaZ} or a more recent source such as \cite{Gut} for a detailed exposition showing how $\rho$ may be continuously extended from the densely defined set of all semisimple $A \in Sp(2n)$, it's purpose here is to guarantee the $E^i_g$ and $E^i_s$ do not become singular in the limit. 
	\end{Remark}
With these three ingredients and the fact that $\rho$ is multiplicative with respect to direct sums we may prove the theorem by defining $\hat{\rho}(L) := \rho^2(\phi)$ (details found in section~\ref{subsec:finalproof}).

\section{Iterating Linear Canonical Relations}
As mentioned above the composition map is not continuous which motivates the following classic example demonstrating the need for transversality.
\begin{Example}\label{example:cont}
	Let $K_i = Gr(A_i)$ and $K_i^\prime = Gr(A_i^{-1})$ for $\left\{A_i\right\}_{i=1}^\infty \subset Sp(2n)$ such that $K_i \to K = L_1 \times \left\{0\right\} \oplus \left\{0\right\} \times L_2  \in \Lambda_{2n}$ where the $L_i \in \Lambda_n$ for $i=1,2$. Then the set-theoretic composition $K_i^\prime \circ  K_i  = \triangle_V$ for all $i \in \mathbb{N}$ so that $lim_{i\to \infty}(K_i^\prime \circ  K_i) = \triangle_V$. Yet we observe $K_i^\prime \to K^\prime = L_2 \times \left\{0\right\} \oplus \left\{0\right\} \times L_1$ as $i \to \infty$ so that the composition of their limits shows a failure of continuity, 
	\begin{eqnarray}
	lim_{i\to \infty} K_i^\prime \circ lim_{i\to\infty} K_i =  L_1 \times \left\{0\right\} \oplus \left\{0\right\} \times L_1 \neq \triangle_V = lim_{i\to \infty}(K_i^\prime \circ  K_i)\nonumber.
	\end{eqnarray}
	
	Note here that $ran(K) = dom(K^\prime)$ so the pair are in some sense maximally non-transversal.
\end{Example}

Regarding $L \in \Lambda_{2n}\setminus H$ we have $dom(L) \oplus halo(L) = V$ so certainly $dom(L) + ran(L) = V$ and $dom(L)$ is transversal to $ran(L)$. This suggests the construction of a map for any $i \in \mathbb{Z}$, 
\begin{Definition}\label{def:sqr} For any $i \in \mathbb{Z}$ define composition as follow,
	\begin{eqnarray}
	(\ast)^i\colon \Lambda_{2n}\setminus H &\to& \Lambda_{2n}\setminus H\nonumber\\
	L &\mapsto& L^i = \underset{i \text{ times}}{\underbrace{L \circ L \circ \dots \circ L}}.\nonumber 
	\end{eqnarray}
	where $L^0 := \triangle_{V} = \left\{(v,v) \in V \times V \, | \, v\in V\right\}$. 
\end{Definition}

The case when $k = 0$ is justified since $\triangle_V$ is an `identity' for linear Lagrangian relations; that is $\triangle_V \circ L = L = L \circ \triangle_V$ \cite{Wei}.

\begin{Lemma}\label{lem:sqr}
	Given $L \in \Lambda_{2n}\setminus H$ then $L^i \in \Lambda_{2n}\setminus H$ for all $i \in \mathbb{Z}$ so that $(\ast)^i\colon \Lambda_{2n}\setminus H \to \Lambda_{2n}\setminus H$ is well defined. More specifically for any $L \in \Lambda_{2n} \setminus H$ when $i>0$ we have $dom(L^i) = dom(L)$ and $range(L^i) = range(L)$ and when $i < 0$ we have $dom(L^i) = Ran(L)$ and $range(L^i) = dom(L)$.
\end{Lemma} 

\textit{Proof: }To show this we first observe that $ker(L) \leq dom(L^i)$ for any $i \ge 1$ since $(0,0)$ is contained in every canonical relation. Then if $dim(ker(L)) = k$ we may write a basis $(d_1,\dots,d_{2n-2k})$ such that $\langle d_j\rangle_{j=1}^{2n-2k} \oplus ker(L) = dom(L)$ and thus each $d_j$ is associated (non-uniquely) via $L$ to some $r_j \in range(L)$. The $r_j$ are also linearly independent since if $r_1 = \sum_{j=2}^{2n-2k} c_jr_j$ then $d_1 - \sum_{j=2}^{2n-2k} c_jd_j \in ker(L)$ which violates $\langle d_j\rangle_{j=1}^{2n-2k} \pitchfork ker(L)$. 

Now again since $L\not\in H$ then $\langle r_j\rangle_{j=1}^{2n-2k} + halo(L) = ran(L)$ and we see that $dom(L) \cap ran(L) = dom(L) \cap \langle r_j\rangle_{j=1}^{2n-2k} $. Since $dim(dom(L) \cap ran(L)) \ge 2n-2k$ we see for dimensional reasons that $\langle r_j\rangle_{j=1}^{2n-2k} \leq dom(L)$ and therefore 
\begin{eqnarray}
dom(L^2) = (\langle r_j\rangle_{j=1}^{2n-2k}) \oplus ker(L) = dom(L).\nonumber
\end{eqnarray} 

An identical argument shows that $ran(L^2) = ran(L)$ as well. To conclude for $i \in \mathbb{N}$ we use this as the base case of a simple inductive argument regarding the domain and range of $L^i \circ L$ and $L \circ L^i$ which suffices to prove the claim for $i \ge 2$.

As for negative powers we set $L^{-1} := L^r$ where $(v,w) \in L^r \Leftrightarrow (w,v) \in L$, motivated by the fact that when $L = Gr(\phi)$ then $L^{r} = Gr(\phi^{-1})$. Since the invariants described in the next section are symmetric with respect to the domain and range it follows that $\kappa(L^r) = \kappa(L)$ (see ~\eqref{eqn:invariants}) implying $L \in \Lambda_{2n}\setminus H \Leftrightarrow L^r \in \Lambda_{2n} \setminus H$ so that $L^{-i} := (L^r)^i$ is well defined for any $i \in \mathbb{N}$ concluding the proof for all $i \in \mathbb{Z}$ given the special definition when $i = 0$. $\qed$

We prove in lemma~\ref{lem:rhohomo} that this map is continuous for each $i \in \mathbb{N}$. 

\begin{Remark}
extending homogeneity of the extended mean index over paths to negative numbers becomes a bit more interesting as it reverses the isotropic pair of $L$ covered in the next section, steering away from this we will assume the powers to be non-negative as stated in the theorem for the extended mean index.	
\end{Remark} 

\section{Isotropic Pairs}
\subsection{Conjugacy Classes of Isotropic Pairs}\label{subsub:iso}

Denote the Grassmannian of isotropic subspaces of $(V,\omega)$ with dimension $k$ by $I_{k}(V)$ and call any ordered pair $(B_1,B_2) \in I_{k}(V)\times I_{k}(V)$ an isotropic pair and consider the following notion of equivalence.
\begin{Definition}
	$(B_1,B_2) \sim (B^\prime_1,B^\prime_2)$ if and only if there exists $A \in Sp(V)$ for which $(A(B_1),A(B_2)) = (B^\prime_1,B^\prime_2)$. The equivalence classes coincide with the orbits of the group action $Sp(V)\circlearrowright I_k(V) \times I_k(V)$ where $A \cdot (B_1,B_2) = (A(B_1),A(B_2))$.
\end{Definition}
For our purposes we have assumed the isotropic pairs have the same dimension; our goal in introducing them is to examine $L \in \Lambda_{2n}$ via the associated isotropic pair $(B_1,B_2) = (ker(L),halo(L))$ and since $dim(ker(L)) = dim(halo(L))$ this assumption is justified.
As shown in \cite{LorW} the following four integers (with the dimensional constraint) form a complete set of invariants for isotropic pairs subject to the relations $0 \leq r \leq \kappa \leq k \leq n$ and $0 \leq \kappa-r \leq n-k$, 
\begin{eqnarray}\label{eqn:invariants}
(\kappa,r,k,n) = \left(dim(B_1^\omega \cap B_2),dim(B_1 \cap B_2), dim(B_{1}), \frac{1}{2}dim(V)\right).
\end{eqnarray}

Now denoting $\Lambda^k_{2n} := \left\{L \in \Lambda_{2n} \, | \, dim(ker(L)) = k\right\}$ for any $0 \leq k \leq n$ and noting $H = \left\{ L \in \Lambda_{2n} \, | \, \kappa(L) \ge 1\right\}$ we see that the above equivalence relation on $I_k(V)\times I_k(V)$ induces an equivalence relation on $\Lambda_{2n}^k$ for each $k \leq n$ (and therefore on all of $\Lambda_{2n}$) where $L \sim L^\prime$ if and only if $(ker(L),halo(L)) \sim (ker(L^\prime),halo(L^\prime))$. A detail to note is that on $\Lambda_{2n}^0 \cong Sp(V)$ all maps belong to a single equivalence class under this equivalence relation. We may compare this equivalence relation induced by the $Sp(2n)$ action on isotropic pairs to a finer relation on $\Lambda_{2n}$ induced by an essentially identical $Sp(2n)$ action now acting on $\Lambda_{2n}$.
\begin{Definition}
	$L \sim_{Gr} L^\prime$ if and only if $(v,w) \in L \Leftrightarrow (Av,Aw) \in L^\prime$. The equivalence classes coincide with the orbits of the group action $Sp(V)\circlearrowright \Lambda_{2n}$ where $A\cdot((x,y) \in L) \mapsto (Ax,Ay) \in A\cdot L$.
\end{Definition}
This equivalence relation in particular splits the single $\sim$ equivalence class of $\Lambda_{2n}^0$ into the usual conjugacy classes of the underlying symplectomorphism for each $L \in\Lambda^0_{2n}$ while conversely $\sim$ and $\sim_{Gr}$ are identical on $\Lambda_{2n}^n$. The classification and production of normal forms for $L \in \Lambda_{2n}$ with respect to this finer equivalence relation is, to the author's knowledge, close to completed \cite{Wei18} following the partial results of \cite{Lor}.

\subsection{The Codimension of $H$}

First we will need a lemma found in \cite{He},

\begin{Lemma}\label{lem:Fibration}
	Given any $L \in \Lambda_{2n}^k$ there exists a canonical fibration:
	\begin{eqnarray}
	\Lambda_{2n-2k}^0(dom(L)/ker(L) \times ran(L)/halo(L)) \hookrightarrow \Lambda_{2n}^k(V \times \overbar{V}) &\underset{Pr}{\twoheadrightarrow}& I_k(V) \times I_k(V)\nonumber\\
	L &\mapsto& (ker(L),halo(L))\nonumber
	\end{eqnarray}
\end{Lemma}
\begin{Remark}
	From this result we can define 
	\begin{eqnarray}
	\hat{H}_k := \left\{(B_1,B_2) \in I_k(V) \times I_k(V) \, | \, dim(B_1^\omega\cap B_2) \ge 1\right\}\nonumber
	\end{eqnarray} so that $H \cap \Lambda^k_{2n} = Pr^{-1}(\hat{H}_k)$ for each $1 \leq k\leq n-1$. When $k = n$ the fibers of the above fibration collapse leaving $\hat{H}_n \subset \Lambda^n_{2n} \cong I_n(V) \times I_n(V)$. In particular each $H \cap \Lambda^k_{2n}$ is a principal $Sp(2n-2k)$ bundle over $\hat{H}_k$ for $1 \leq k \leq n-1$.
\end{Remark}
\begin{Remark}
	Canonical in this case refers to the fact that the bundle components are all derived from a shared symplectic vector space so that the fibers are smoothly dependent on the base point. The total space $\Lambda_{2n}^k(V\times \overbar{V})$ is a stratum of the Lagrangian Grassmannian yielding the fixed symplectic vector space $(V\times \overbar{V},\tilde{\omega})$. The base is a product of two isotropic Grassmannians over $(V,\omega)$ such that the following symplectic vector space is determined uniquely by the base point $(B_1,B_2) \in I_k(V)\times I_k(V)$,
	\begin{eqnarray}
	\left(B_1^\omega/B_1 \times B_2^\omega/B_2, \tilde{\omega}_{red} = pr^\ast_1\omega_{1,red} - pr^\ast_2\omega_{2,red}\right)\nonumber
	\end{eqnarray}
	over which the fiber is defined, namely $\Lambda_{2n-2k}^0(B_1^\omega/B_1 \times B_2^\omega/B_2) \cong Sp(2n-2k)$. So then since $\tilde{\omega}_{red}$ is the co-isotropic reduction of $\left.\tilde{\omega}\right|_{B_1^\omega \times B_2^\omega}$ for fixed $k$ the vector spaces and their induced symplectic form over which the fibers are defined vary smoothly over the base space as subspaces of the ambient vector space.
\end{Remark}

\begin{Theorem}\label{thm:codim}
	$codim(H) = 2$.
\end{Theorem}

\textbf{Proof:} We refer to \cite{He} where it is shown that $codim(\Lambda_{2n}^k) = k^2$ so that since $H \cap \Lambda_{2n}^0 = \emptyset$ we have the trivial bound of $1 \leq codim(H) $. To sharpen this we consider the equivalence classes $[L] \subset \Lambda_{2n}^k$ induced by the associated isotropic pair class $[(ker(L),halo(L))] \subset I_k(V) \times I_k(V)$ when $k=1$. There are three equivalence classes contained in $I_1(V) \times I_1(V)$ (and therefore three in $\Lambda_{2n}^1$) but the class $(\kappa,r,k) = (0,0,1)$ does not intersect $H$ so that only the two classes satisfying $\kappa = k = 1$ need be checked, namely $dim(ker(L) \cap halo(L)) = r = 0$ or $r=1$.

\textbf{Case I: r=1} 

Let $L \in [L]_1 = \left\{ L \in \Lambda^1_{2n} \, | \, \kappa(L) = r(L) = 1 \right\}$. With the above lemma we first consider the associated class $[(ker(L),halo(L))]$ and write $dom(L) = halo(L) =\langle v \rangle$ for any $v \in V$. Then since $v$ is arbitrary and all one dimensional subspaces are isotropic we see that $[(\langle v\rangle,\langle v\rangle)] = \triangle_{I_1(V)} \subset I_1(V) \times I_1(V) \cong \mathbb{RP}^{2n-1} \times \mathbb{RP}^{2n-1} $ implying that $dim([(\langle v\rangle,\langle v\rangle)]) = 2n-1$. Thus it follows that 
\begin{eqnarray}
dim([L]_1) &=& dim(Sp(2n-2)) + dim([(\langle v\rangle,\langle v\rangle)])\nonumber \\
&=& 2n^2 - 3n + 1 + (2n -1)\nonumber\\
&=& 2n^2 - n = dim(\Lambda_{2n}) - 2n\nonumber
\end{eqnarray}
so that $codim([L]_1) = 2n$ in $\Lambda_{2n}$.
\newline
\newline
\newline
\textbf{Case II: r=0}

When $L \in [L]_0 = \left\{L \in \Lambda_{2n}^1 \, | \, r(L) = 0, \kappa(L) = 1\right\}$ we see that 
\begin{eqnarray}
[(ker(L),halo(L))] = \left\{ (B_1,B_2) \, | \, B_1 \neq B_2, \text{   and   } B_2 \leq B_1^\omega\right\}\nonumber
\end{eqnarray}
where the first condition is due to $r(L) = 0$ and the second from $\kappa(L) =1$. As before there are $2n-1$ dimensions in freely choosing $B_1$ while the first and second conditions together imply that $B_2$ is restricted to the image of $B_1^\omega \setminus B_1$ (a hyperplane missing a one dimensional subspace) which descends under the quotient map $\mathbb{R}^{2n} \to \mathbb{RP}^{2n-1}$ to a punctured projective hyperplane $[B_1^\omega]\setminus \left\{[B_1]\right\} \subset \mathbb{RP}^{2n-1}$ so that there are $2n-2$ dimensions available when choosing $B_2$. This yields $dim([(ker(L),halo(L))]) = 4n - 3$ so that
\begin{eqnarray}
dim([L]_0) = dim(Sp(2n-2)) + 4n - 3 = 2n^2 + n - 2 = dim(\Lambda_{2n}) - 2\nonumber.
\end{eqnarray}

Since these are the two equivalence classes in the stratum of minimal codimension intersecting $H$ we see that $codim(H) = 2$. $\qed$

\section{Decomposing Linear Canonical Relations in the Complement of $H$}

Given some $L \in \Lambda_{2n}$ we have the example~\ref{example:squared} showing that $\rho$ may not be continuously extended to $\Lambda_{2n}\setminus H$. Following section~\ref{subsub:iso} we see for $1 \leq k \leq n$ that the $L \in \Lambda_{2n}^k$ such that $dom(L) \pitchfork halo(L)$ form a dense set in $\Lambda_{2n}^k$ and as shown in \cite{G15} this condition yields a natural domain on which $\rho^2$ may be extended continuously. Conversely in example~\ref{example:converse} we produce an $L \in H$ for which there exists a uniquely defined graph component. The conjugacy classes described in theorem~\ref{thm:dimbound} on which $L$ possesses a graph component with a unique symplectic map (up to a choice of domain) suggests that $\Lambda_{2n}\setminus H$ is potentially not the maximal subset on which $\hat{\rho}$ may be extended continuously.

In the following section we prove that the $L \not\in H$ induce a unique $\omega$-orthogonal decomposition of $V$ which informally splits $L$ into the direct sum of its `singular' and `graph' components.

\textit{Note: We will proceed denoting $L_1 := dom(L)$ and $L_2 := ran(L)$ for a given $L \in \Lambda_{2n}$.}

\begin{Theorem}\label{subthm:decomp}
	Given $L \in \Lambda_{2n}$ satisfying $L_1 \cap L_2^\omega = \left\{0\right\}$ there exists a unique symplectic decomposition $V = V_s \oplus V_g$ and symplectic map $\phi\colon V_g\circlearrowleft$ for which 
	\begin{eqnarray}
	L = (L_1^\omega \times \left\{0\right\} \oplus \left\{0\right\} \times L_2^\omega)\oplus Gr(\phi) \leq (V_s \times \overbar{V_s}) \oplus (V_g \times \overbar{V_g})\nonumber
	\end{eqnarray}
	where the $L_i^\omega$ are transverse Lagrangian subspaces of $V_s$.
	
\end{Theorem}

\textbf{Proof:} 

Given every $(v,w) \in L$ we may construct a map: 
\begin{eqnarray}
\tilde{\phi}\colon L_1 &\to& L_2/L_2^\omega\nonumber\\
v &\mapsto& [w]\nonumber.
\end{eqnarray}

\textbf{(i)} This map is well defined. \textit{Proof:} Given $(v,w),(v,w^\prime) \in L$ we have that $[w] = [w^\prime] \Leftrightarrow w-w^\prime \in L_2^\omega$. 

\textbf{(ii)} The kernel of $\tilde{\phi}$ is $L_1^\omega$. \textit{Proof:} Since $\tilde{\phi}(v) = [w] = 0 \Leftrightarrow (v,w) \in L, \,\, w \in L_2^\omega \Leftrightarrow (0,w) \in L$. Then by linearity we see $(v,0) \in L$ which by definition means $v \in L_1^\omega$. Conversely if $v \in L_1^\omega$ then $(v,0) \in L \Rightarrow \tilde{\phi}(v) = [0]$ so indeed $ker(\phi) = L_1^\omega$.

Since the $L_i$ are co-isotropic the map $\phi\colon L_1/L_1^\omega \to L_2/L_2^\omega$ is an isomorphism between symplectic vector spaces where for both $i=1,2$ we have that
\begin{eqnarray}
dim(L_i/L_i^\omega) = dim(L_i) - dim(L_i^\omega) = (2n-k)-k = 2n-2k. \nonumber 
\end{eqnarray}

In fact since the $L_i/L_i^\omega$ are reduced co-isotropic subspaces they each possess a canonical symplectic form: $\omega^i_{red}([v],[v^\prime]) := \omega(v,v^\prime)$ for all $v,v^\prime \in L_i$ which is independent of the choice of representatives $v,v^\prime \in L_i$ for $i=1,2$. This leads us to our next step.

\textbf{(iii)} The map $\phi\colon (L_1/L_1^\omega,\omega^1_{red}) \to (L_2/L_2^\omega,\omega^2_{red})$ is symplectic. \textit{Proof:} Given any pair $(v,w),(v^\prime,w^\prime) \in L$ they must satisfy 
\begin{eqnarray}
\tilde{\omega}((v,w),(v^\prime,w^\prime)) = 0 \Leftrightarrow \omega(v,v^\prime) = \omega(w,w^\prime).\nonumber
\end{eqnarray} It follows for any $[v], [v^\prime] \in L_1/L_1^\omega$ and $[w],[w^\prime] \in L_2/L_2^\omega$ such that $(v,w),(v^\prime,w^\prime) \in L$ that $\omega_{red}([v],[v^\prime]) = \omega_{red}([w],[w^\prime]) = \omega_{red}(\phi[v],\phi[v^\prime])$ so $\phi$ is indeed a symplectic map between the two reduced spaces. 

There always exists a pair of symplectic subspaces $V_i \leq L_i$ which are mapped bijectively under the projection maps $\pi_i\colon L_i \to L_i^\omega$ so that $V_i \cong L_i/L_i^\omega$ and $\phi\colon V_1 \to V_2$ satisfies $\phi^\ast \omega = \omega$. This map, although symplectic depends not only on $L$ but which \textit{pair} of $V_i$ are chosen as well. 

Until this point the hypothesis that $\kappa(L) = 0$ has not been needed but is now required to produce a unique $V_g := V_1 = V_2$ on which $L$ induces $\phi \in Sp(V_g) \cong Sp(2n-2k)$.

\begin{Sublemma}\label{lem:graph}
	If $\kappa(L) = 0$ then $V_g := (L_1 \cap L_2,\omega)$ is a symplectic subspace of $V$ and $(V_g,\omega) \cong (L_i/L_i^\omega,\omega^i_{red})$ for $i = 1,2$.
\end{Sublemma}

\textbf{Proof:} First note the following three conditions are sufficient to show the above: \textbf{(1)} $dim(L_1 \cap L_2) = 2n-2k$, \textbf{(2)} $V_g$ is a symplectic subspace of $V$ and \textbf{(3)} $ker(\pi_i) = L_i^\omega$ has trivial intersection with $V_g$ for $i = 1,2$. if \textbf{(1)}-\textbf{(3)} are satisfied then $\left.\pi_i\right|_{V_g}$ is a symplectomorphism for both $i = 1,2$.

\textbf{(1)} $dim(L_1 \cap L_2) = 2n-2k$
	
	\textit{Proof:} Recalling that $dim(L_i) = 2n -k$ for some $k \leq n$ meaning $dim(L_i^\omega) = k$ and we see from $L_1^\omega \leq L_1$ that
	\begin{eqnarray}
	\qquad\qquad\kappa(L) = dim(L_1\cap L_2^\omega) = 0 \Rightarrow dim(L_1^\omega \cap L_2^\omega) = 0 \Leftrightarrow dim(L_1^\omega \oplus L_2^\omega) = 2k\nonumber
	\end{eqnarray} 
	and $dim(L_1\cap L_2) = 2n - dim(L_1^\omega \oplus L_2^\omega) = 2n-2k$.
	
\textbf{(2)} $V_g$ is a symplectic subspace of $V$
	
	\textit{Proof:} We note $\kappa(L) = dim(L_1 \cap L_2^\omega) = dim(L_1^\omega \cap L_2) = 0$ implies via inclusion that
	\begin{itemize}
		\item $L_1\cap L_2^\omega = \left\{0\right\} \Rightarrow (L_1 \cap L_2) \cap L_2^\omega = \left\{0\right\}$ \\
		\item $L_2\cap L_1^\omega = \left\{0\right\} \Rightarrow (L_1 \cap L_2) \cap L_1^\omega = \left\{0\right\}$
	\end{itemize}
	
	so that since $2n-2k + 2k = 2n$ we have	the following decomposition,
	\begin{eqnarray}
	V = V_s \oplus V_g := (L_1^\omega \oplus L_2^\omega) \oplus (L_1 \cap L_2).\nonumber
	\end{eqnarray}
	
	Indeed since $V_g^\omega = V_s$ and the above spans $V$ we see that the two form a pair of complementary symplectic subspaces depending uniquely on $L$ (or more precisely the isotropic pair associated to $L$).

\textbf{(3)} $V_g \cap L_i^\omega = \left\{0\right\}$ for $i = 1,2$. \textit{Proof:} This is an immediate consequence of the above decomposition. $\qed$

Thus letting $V_g := L_1 \cap L_2$ we see since $(\pi_i^{\ast}\left.\omega_{red})\right|_{V_g \times V_g} = \left.\omega\right|_{V_g \times V_g}$ for $i = 1,2$ there is a $\phi \in Sp(V_g)$ such that $Gr(\phi) \leq V_g \times \overbar{V_g}$ is Lagrangian. Then since the $L_i^\omega$ are isotropic in $V$ of dimension $k$ it follows that the $L_i^\omega$ are transverse maximal isotropic  subspaces of $V_s$ proving the last claim of theorem~\ref{subthm:decomp}. $\qed$

\section{The Circle Map $\rho$}\label{sec:rho}

\subsection{Properties of $\rho$}

In addition to the mean index the circle function $\rho$ has been used (as it was in \cite{SaZ}) to construct the Conley-Zehnder index.

\begin{Definition}\label{def:rho}\cite{SaZ}
	We define $\rho\colon Sp(2n) \to S^1$ as follows. Given $A \in Sp(2n)$ let $E = Spec(A) \cap (S^1 \cup \mathbb{R})$ be the collection of real and elliptic eigenvalues of $A$. For elliptic eigenvalues $\lambda \in E\cap(S^1 \setminus \left\{\pm 1\right\})$ define $m^+(\lambda)$ to be the number of positive eigenvalues of the symmetric, non-degenerate two form $Q$ defined on the complex eigenspace $E_\lambda$ where
	\begin{eqnarray}
	Q\colon E_\lambda \times E_\lambda &\to& \mathbb{R}\nonumber\\
	(z,z^\prime) &\mapsto& Im(\omega(z,\overbar{z^\prime})).\nonumber
	\end{eqnarray} 
	
	Then letting $m^-$ denote the sum of the algebraic multiplicities for the real negative eigenvalues we have
	\begin{eqnarray}
	\rho(A) := (-1)^{\frac{1}{2}m^-} \prod_{\lambda \in S^{1} \setminus \left\{\pm 1\right\}} \lambda^{\frac{1}{2}m^+(\lambda)}. \nonumber
	\end{eqnarray}
\end{Definition}

In our case since eigenvalues are unique the term $m^+(\lambda) = 1$, the counting of positive eigenvalues of $Q$ over the eigenspace for $\lambda$ amounts to a consistent way of choosing a single value from an elliptic eigenvalue pair while still taking the product over all elliptic eigenvalues. 
\begin{Proposition}\cite{SaZ}
	
	The map $\rho\colon Sp(2n) \to S^1$ has the following properties:
\end{Proposition}
\begin{enumerate}
	\item \textit{(determinant)} For $A \in U(n) \subset Sp(2n)$ we have $\rho(A) = Det_{\mathbb{C}} (A)$.
	\item \textit{(invariance)} $\rho$ is invariant under conjugation, 
	\begin{eqnarray}
	\rho(B^{-1}AB) = \rho(A), \,\, \forall B \in Sp(2n).\nonumber
	\end{eqnarray}
	\item \textit{(normalization)} $\rho(A) = \pm 1$ if $A$ has no elliptic eigenvalues.
	\item \textit{(multiplicativity)} If $A = B\oplus C \in Sp(2n) \times Sp(2m) \subset Sp(2(n+m))$ then $\rho(A \oplus B) = \rho(A) \rho(B)$.
	\item \textit{(homogeneity)} If $A \in Sp(2n)$ we have that $\rho(A^l) = \rho(A)^l$ for any $l \in \mathbb{Z}$.
\end{enumerate}

All of the above properties are inherited by $\hat{\rho}$ when $L\in \Lambda_{2n}^0$ but for non-graph Lagrangian subspaces some properties no longer have an analog. In example~\ref{example:fundamental} we see that $\hat{\rho}$ is indeed not a circle map on $\Lambda_{2}\setminus H$ since there exists a non-contractible loop $\gamma:I \to \Lambda_2 \setminus H$ for which $\hat{\Delta}(\gamma) = 0$.

\subsection{Properties of the Extension of $\rho$}

Despite the above observation, $\hat{\rho}$ does inherit the above properties with the caveat that they are all only defined for $L \in \Lambda_{2n}^0$ and are in that sense trivial. One non-trivial property is found in the following lemma.

\begin{Lemma}\label{lem:rhohomo}
	\qquad
	\begin{enumerate}
		\item The operation $(\ast)^l$ (definition~\ref{def:sqr}) is a continuous map for any $l \in \mathbb{N}$.
		\item (Homogeneity) Given any $L \in \Lambda_{2n}\setminus H$ with graph part $Gr(\phi)$ and assuming theorem~\ref{thm:rhothm} is true, that is $\hat{\rho}$ is continuous and $\hat{\rho}(L) = \rho^2(\phi)$, then $\hat{\rho}(L^l) = \hat{\rho}(L)^l$.
	\end{enumerate} 
\end{Lemma}

\textbf{Proof:} As shown in lemma~\ref{lem:sqr}, the iterated composition operation restricted to $\Lambda_{2n}\setminus H$ is a well defined map and we may compute in coordinates $L^2$ where
\begin{eqnarray}
L = (Ker(L) \times \left\{0\right\} \oplus \left\{0\right\} \times halo(L))\oplus Gr(\phi)\nonumber
\end{eqnarray}
so that we may verify the claim via a Darboux basis adapted to $V = ker(L) \oplus halo(L) \oplus (dom(L) \cap ran(L))$.

Let $\langle v_i\rangle_{i=1}^k = ker(L)$ and $\langle w_i\rangle_{i=1}^k = halo(L)$. Then since $(v_i,0),(0,w_j) \in L$ for all $i,j \leq k$ it follows that $(v_i,w_j) \in L \circ L$ for all $i,j \leq k$, namely $ker(L) \times \left\{0\right\} \oplus \left\{0\right\} \times halo(L) \leq L \circ L$. Next consider some $(v,w) \in Gr(\phi)$ so that since $\phi \in Sp(V_g)$ for any $w \in V_g$ there exists a unique $z \in V_g$ for which $(w,z) \in Gr(\phi)$ implying $(v,z) \in L \circ L$. Intuitively this states that $Gr(\phi) \circ Gr(\phi) = Gr(\phi^2)$. As for $(0,w_i) \in halo(L)$ in the first $L$ and $(v_i,0) \in ker(L)$ belonging to the second $L$ the only resulting vector derived from these in the product is $(0,0)$ (regardless of whether $L \in H$ or not) so we see for any $l \ge 1$ that
\begin{eqnarray}
L^l := (Ker(L) \times \left\{0\right\} \oplus \left\{0\right\} \times halo(L))\oplus Gr(\phi^l)\nonumber
\end{eqnarray}
and $\hat{\rho}(L^l) = \hat{\rho}(L)^l$ since $\rho$ is homogeneous on $Sp(V_g)$.

Continuity follows on each $\Lambda_{2n}^k\setminus H$ after referring to lemma~\ref{lem:Fibration}, that is the fibration $Sp(V_g) \hookrightarrow \Lambda^k_{2n}\setminus H  \twoheadrightarrow I_k(V)\times I_k(V)$ for each $1 < k < n$. Since $L, L^l$ share the same domain and range for all $l \in \mathbb{N}$ then $(\ast)^l$ preserves the fibers when $l \ge 1$ (yet swaps the isotropic pair to its reversal for $l \leq -1$). Since the fibers vary smoothly over the base we see that $(\ast)^l$ is continuous since it is continuous on each fiber, inheriting the group operation on $Sp(V_g)$. When $k=0$ this operation corresponds to the group operation in $Sp(2n)$ and when $k = n$ it is the identity map so $(\ast)^l$ is continuous on $\Lambda_{2n}\setminus H$ . 
$\qed$

\begin{Proposition}
	The map $\hat{\rho}\colon \Lambda_{2n} \setminus H \to S^1$ inherits the following properties:
\end{Proposition}
\begin{enumerate}
	\item \textit{(determinant)} For $A \in U(n) \subset Sp(2n)$ we have that $\hat{\rho}(Gr(A)) = (Det_{\mathbb{C}} (A))^2$.
	\item \textit{(invariance)} $\hat{\rho}$ is invariant under conjugation on $\Lambda_{2n}^0$, 
	\begin{eqnarray}
	\hat{\rho}(Gr(B^{-1}AB)) = \hat{\rho}(Gr(A)) \,\, \forall B \in Sp(2n).\nonumber
	\end{eqnarray}
	\item \textit{(normalization)} $\hat{\rho}(Gr(A)) = 1$ if $A$ has no elliptic eigenvalues.
	\item \textit{(multiplicativity)} If $L = K_1\oplus K_2 \in \Lambda_{2n}\setminus H \times \Lambda_{2m}\setminus H \subset \Lambda_{2(n+m)}\setminus H$ then $\hat{\rho}(K_1\oplus K_2) = \hat{\rho}(K_1)\hat{\rho}(K_2)$.
\end{enumerate}
	\textbf{Proof (4):} We begin by observing when $L = K_1 \oplus K_2 \in \Lambda_{2n}\setminus H \times \Lambda_{2m}\setminus H$ there exists a pair of symplectic bases, $(x_i,y_i)_{i=1}^n$ and $(u_i,v_i)_{i=1}^m$ which together form subspaces which form a decomposition isomorphic to the induced one: $\mathbb{R}^{2n} \oplus \mathbb{R}^{2m}  = \mathbb{R}^{2(n+m)}$. With this decomposition we may write
	\begin{eqnarray}
	L &= (ker(K_1) \times \left\{0\right\} \oplus \left\{0\right\} \times halo(K_1) \oplus Gr(\phi_{K_1}))\nonumber \\
	&\oplus\nonumber \\
	&(ker(K_2) \times \left\{0\right\} \oplus \left\{0\right\} \times halo(K_2) \oplus Gr(\phi_{K_2})).\nonumber
	\end{eqnarray}

	Since $dom(K_1),ran(K_1) \leq \langle x_i,y_i\rangle_{i=1}^{n}$ and $dom(K_2),ran(K_2) \leq \langle u_i,v_i\rangle_{i=1}^m$ then theorem~\ref{subthm:decomp} implies the two subspaces $V^j_s$ and $V^j_g$ determined by each of the $K_1$ and $K_2$ have pairwise trivial intersection thereby refining the decomposition,
	\begin{eqnarray}
	 \mathbb{R}^{2(n+m)} = (V^1_s \oplus V^1_g) \oplus (V^2_s \oplus V^2_g) = \mathbb{R}^{2n} \oplus \mathbb{R}^{2m}. \nonumber
	\end{eqnarray}
	
	This implies $ker(L) = ker(K_1) \oplus ker(K_2)$, $halo(L) = halo(K_1) \oplus halo(K_2)$ and $Gr(\phi_L) = Gr(\phi_{K_1}) \oplus Gr(\phi_{K_2})$, in particular this means that $\phi_L = \phi_{K_1} \oplus \phi_{K_2} \in Sp(2n) \times Sp(2m) \subset Sp(2(n+m))$. To conclude the proof we assume theorem~\ref{thm:rhothm} is true so that $\hat{\rho}$ is multiplicative since $\hat{\rho}(L) = \rho^2(\phi_L) = \rho^2(\phi_{K_1})\rho^2(\phi_{K_2}) = \hat{\rho}(K_1)\hat{\rho}(K_2)$.$\qed$

\section{The Asymptotic Behavior of Unbounded Sequences in $Sp(2n)$}
\subsection{A Sufficient Condition for Asymptotic Hyperbolicity}

The following theorem states that any $A \in Sp(2n)$ with $Gr(A)$ sufficiently near $L = L_1\times \left\{0\right\} \oplus \left\{0\right\} \times L_2$ with $L_1 \pitchfork L_2$ has only hyperbolic eigenvalues.

\begin{Theorem}\label{subthm:hyperbolic}
	Suppose $\left\{A_i\right\}_{i=1}^\infty \subset Sp(2n)$ is a sequence of symplectomorphisms each with distinct eigenvalues such that $Gr(A_i) \underset{i \to \infty}{\longrightarrow} L_1\times \left\{0\right\} \oplus \left\{0\right\} \times L_2 \in \Lambda_{2n}$, $L_i \in \Lambda_n$ and $L_1 \cap L_2 = \left\{0\right\}$. Then there exists $K\in \mathbb{N}$ such that 
	\begin{eqnarray}
	Spec(A_i) \cap (S^1\setminus \left\{\pm 1\right\}) = \emptyset \nonumber
	\end{eqnarray} for all $i \ge K$.
\end{Theorem}	

\textbf{Proof:} Assume there exists some $\lambda_i \in Spec(A_i)$ such that $\lambda_i \in  S^1\setminus \left\{\pm 1\right\}$ for all $i \ge K \in \mathbb{N}$. Then since there are no multiple roots such a $\lambda_i$ belongs to a unique symplectic eigenvalue pair $\left\{\lambda_i,\overbar{\lambda_i}\right\}$ with $|\lambda_i| = 1$ which are stable in the sense that they remain elliptic away from the exceptional points $\pm 1$ since the points $\pm 1$ are the only values at which a (unique) elliptic eigenvalue pair may become hyperbolic or in general meet another pair to form a quadruple upon passing to the limit (this is certainly not true without uniqueness). For more details on how to go about `ordering' the eigenvalues of a sequence of maps so that individual sequences of eigenvalues may be coherently formed (as done above implicitly) refer to section~\ref{subsec:convergence}.

This pair of eigenvalues has eigenvectors $x_i\pm iy_i \in \mathbb{C}^{2n}$ with a convergent subsequence of corresponding real eigenspaces $E_i = \langle x_i,y_i\rangle$ with $E_i \to  E$ as $i \to \infty$, on which $A_i$ is conjugate to a rotation for all $i \in \mathbb{N}$. We obtain individual limit vectors by letting $x := \lim_{i \to \infty}\frac{x_i}{|x_i|}$ and $y = \lim_{i \to \infty}\frac{y_i}{|y_i|}$ so we may write $\langle x,y \rangle = E$.

The $E_i$ are symplectic so there exists a decomposition $V = E_i \oplus F_i$ where the $F_i$ form a sequence of some symplectic complements to each $E_i$ and $A_i = \psi_i \oplus \phi_i\colon E_i \oplus F_i \to E_i \oplus F_i$ where $\psi_i$ and $\phi_i$ are symplectic for each $i \in \mathbb{N}$. 

\begin{Sublemma}{The Limit Lagrangian's Kernel and Halo.}\label{lem:asymp} 
	
	Let $\left\{A_i\right\}_{i=1}^\infty \subset Sp(2n)$ denote a sequence of symplectic maps for the following:
	
	\begin{itemize}
		\item For $\left\{A_i\right\}_{i=1}^\infty$ such that $Gr(A_i) \to L \in \Lambda_{2n}\setminus\Lambda_{2n}^0$ then $v \in ker(L) \Leftrightarrow A_iv \to 0$ and $v \in halo(L) \Leftrightarrow A_i^{-1}v \to 0$ as $i \to \infty$.		
		\item For $\left\{A_i\right\}_{i=1}^\infty$ such that such that $Gr(A_i) \to L = L_1 \times \left\{0\right\} \oplus \left\{0\right\} \times L_2$ it is true that $L_1 \cap A^{-1}_iL_2 = \left\{0\right\}$ for sufficiently large $i$.
		
	\end{itemize}
\textit{Note that these both hold regardless of whether $L \in H$ or not.} 
\end{Sublemma}	
\textbf{Proof:} We prove the first claim after which the lemma follows via contradiction:

\begin{enumerate}
	\item $v \in ker(L) \Leftrightarrow A_iv \to 0$. 
	
	\textit{Proof:} We observe $v \in ker(L)$ if and only if $(v,0) \in L$ so that since $Gr(A_i) \to L$ it must be that $(v,A_iv) \to (v,0)$.\\
	\item $v \in halo(L) \Leftrightarrow A_i^{-1}v \to 0$.
	
	\textit{Proof:} Again $v \in halo(L)$ if and only if $(0,v) \in L$. Yet $(A_i^{-1}v,v) \in Gr(A_i)$ so as above we see that $v \in halo(L) \Leftrightarrow A_i^{-1}v \to 0$.
\end{enumerate}

We now prove the second part of the lemma: given some sequence $\left\{A_i\right\}_{i =1}^\infty \subset Sp(2n)$ for which $Gr(A_i) \to L \in \Lambda_{2n}^n$ suppose there exists a sequence $\left\{v_i\right\}_{i=1}^\infty \subset V$ and some $N \in \mathbb{N}$ for which $v_i \to v \neq 0$ and $v_i \in L_1 \cap A_i^{-1}L_2$ for all $i \ge N$. 

Then for any $i \ge N$ we have $v_i \in L_1$ implies $|A_iv_i| \to 0$ and $v_i \in A_i^{-1}L_2$ implies that $w_i := \frac{A_iv_i}{|A_iv_i|} \in L_2$ and therefore $|A_i^{-1}w_i|\to 0$. This means that
\begin{eqnarray}
|A_i^{-1}w_i| = \frac{|v_i|}{|A_iv_i|} \to 0 \nonumber
\end{eqnarray} and since $|v_i| \to |v| \neq 0$ we have the contradiction $|A_iv_i| \to \infty$. $\qed$

Equivalently for large $i$ we see the $A_i$ induce a sequence $L_1 \oplus A^{-1}_iL_2 = V$ of Lagrangian splittings so there exists a unique decomposition for any sequence $u_i = v_{i} + w_{i} \in E_i$ with $v_i \in L_1$ and $w_i \in A^{-1}_iL_2$ such that $|A_iv_i| \to 0$ and $|A_iw_i| \to \infty$. Denote $\psi_i := \left.A_i\right|_{E_i}$ so that $|\psi_i v_i| \to 0$ and $|\psi_i w_{i}| \to \infty$ and consider the following lemma. 
\begin{Sublemma}\label{lem:halo}
	Given a sequence of elliptic eigenspaces $E_i \to E$ and symplectic maps $\left\{\psi_i\right\}_{i=1}^\infty$ as above then $E\cap L_1 \neq \left\{0\right\}$.
\end{Sublemma}
\textbf{Proof:} 
First convergence is a consequence of the fact that $dim(E_i)$ is constant for all $i$ and therefore converges to $E$ along some subsequence \cite{Ka}. Suppose $E \cap L_1 = \left\{0\right\}$, that is every sequence $u_i = v_i + w_i \in E_i$ has $w_i \not\to 0$. We would have for every $u_i= v_{i} + w_{i} \in E_i$ that $|\psi_iu_i|  = |\psi_iv_{i}+\psi_iw_{i}| = |\psi_iw_{i} - (-\psi_iv_{i})| \ge \left||\psi_iw_i| - |\psi_iv_i|\right| \to \infty$ since $|\psi_iv_i| \to 0$ and therefore $|\psi_iu_i| \to \infty$ for any sequence $u_i \in E_i$ with $w_i \not\to 0$.

Now $dim(E_i) = 2$ so that $\omega_i := \left.\omega\right|_{E_i\times E_i}$ is an area form on $E_i$ for each $i \in \mathbb{N}$ so we may choose some sequence of balanced neighborhoods $U_i \subset E_i$ about zero on which 
\begin{eqnarray}
\int_{U_i} \omega_i = 1.\nonumber
\end{eqnarray}
\textit{Note: These $U_i$ may grow without bound but it is of no consequence since we will not pass to the limit.}

Then we have for any $M > 1$ a $K$ for which any normalized sequence $u_i \in E_i$, $u_i \to u \neq 0$ has $|\psi_iu_i| > M$ when $i \ge K$. This implies in particular that for each $M > 1$ there exists a $K \in \mathbb{N}$ for which $M U_i \subset \psi_i U_i$ for every $i \ge K$. It follows then that 
\begin{eqnarray}
1=\int_{U_i}\omega_i< \int_{M U_i}\omega_i \leq \int_{\psi_iU_i} \omega_i \nonumber
\end{eqnarray}
for every $i \ge K$.
Each $\psi_i$ is a symplectomorphism so we have a contradiction by choosing any $M > 1$ so that some $K$ exists for which
\begin{eqnarray}\label{eqn:squeeze}
1 =\int_{U_K} \omega_K=\int_{\psi_K U_K} \omega_K \ge \int_{MU_K} \omega_K > 1.
\end{eqnarray} $\qed$

This lemma may be used in a quick proof of the following.
\begin{Sublemma}
	Given the eigenspaces $E_i \to E$ and $\left\{\psi_i\right\}$ as above then $E\cap L_2 \neq \left\{0\right\}$.
\end{Sublemma}
\textbf{Proof:} We have already established there exists a sequence $v_i \to v$ for which $A_iv_i \to 0$ so we consider any $u_i = v_i + w_i \in E_i \cap (L_1 \oplus A_i^{-1} L_2)$ with $w_i \not\to 0$. Then $A_iu_i = A_iv_i + A_iw_i$ so that each $E_i$ is $A_i$ invariant and each $A_iw_i \in L_2 \cap E_i$. Then since $|A_iv_i| \to 0$ as $i \to \infty$ we see $\lim_{i \to \infty} \frac{A_iu_i}{|A_iw_i|} = \lim_{i \to \infty}\frac{A_iw_i}{|A_iw_i|} \to w \in halo(L) \cap E$. $\qed$

The above lemma and corollary show there exists a sequence of bases $\langle v_i,w_i\rangle \in E_i$ on each real elliptic eigenspace such that $v_i \to v \in L_1$ and $w_i \to w \in L_2$. 
\begin{Sublemma}
	Let $E_i$ be a two dimensional real eigenspace for a complex eigenvalue $\lambda \in S^1 \setminus \left\{\pm 1\right\}$ which exists for sufficiently large $i \in \mathbb{N}$. Then given any sequence $\left\{v_i\right\}$ with each $v_i \in E_i$ and $\frac{v_i}{|v_i|} \to v \neq 0$ we claim,
	\begin{eqnarray}
	v \in ker(L) \Leftrightarrow v \in halo(L).\nonumber
	\end{eqnarray}
\end{Sublemma}	
\textbf{Proof: } Consider $\sigma_i \in Sp(E_i)$ such that $\psi_i = \sigma^{-1}_i \circ R(\theta_i) \circ \sigma_i$ where $R(\theta_i)$ denotes the rotation of the $(x_i^\prime,y_i^\prime) = (\sigma_ix_i,\sigma_iy_i)$ plane $E_i$ by $\lambda_i = e^{i\theta_i}$ so that we may compute for any sequence $\left\{v_i\right\}_{i=1}^\infty$,
\begin{eqnarray}
|\psi^{-1}_iv_i|=|(\sigma_i^{-1} \circ R(\theta_i) \circ \sigma_i)^{-1}(v_i)| &=& |(\sigma_i^{-1} \circ R(\theta_i)^{-1}\circ \sigma_i)(v_i)|\nonumber \\
&=& |(\sigma_i^{-1} \circ R(-\theta_i) \circ \sigma_i)(v_i)|\nonumber 
\end{eqnarray}
and we see $\psi_i^{-1}$ is simply the opposite rotation of $\psi_i$ conjugated by the same matrix $\sigma_i \in Sp(E_i)$. There exists $a_i,b_i \in \mathbb{R}$ such that $v_i = a_i x_i + b_i y_i$ so if we let $v_i^\prime = \sigma_i(v_i)$ then $v_i^\prime := \sigma_i v_i = a_i x_i^\prime + b_i y_i^\prime$ and thus
\begin{eqnarray}
\psi_iv_i &=&(\sigma_i^{-1} \circ R(\theta_i) \circ \sigma_i)(v_i) \nonumber\\
&=&  (\sigma_i^{-1} \circ R(\theta_i))(v_i^\prime) \nonumber\\
&=& (\sigma_i^{-1} \circ R(\theta_i))(a_i x_i^\prime + b_i y_i^\prime) \nonumber\\
&=& \sigma^{-1}_i\left(a_i\cos(\theta_i)x^\prime_i + a_i\sin(\theta_i)y^\prime_i\right) + \sigma^{-1}_i\left(b_i\cos(\theta_i)y^\prime_i - b_i\sin(\theta_i)x^\prime_i\right)\nonumber\\
&=& \left(a_i\cos(\theta_i) - b_i\sin(\theta_i)\right)x_i + \left(a_i\sin(\theta_i) + b_i \cos(\theta_i)\right)y_i\nonumber
\end{eqnarray}
and with squared norm
\begin{eqnarray}
|\psi_iv_i|^2 &=& \left(a_i^2\cos^2(\theta_i) - 2a_ib_i\sin(\theta_i)\cos(\theta_i) + b_i^2\sin^2(\theta_i)\right)\nonumber \\&&+ \left(b_i^2\cos^2(\theta_i) + 2a_ib_i\sin(\theta_i)\cos(\theta_i) +  a_i^2\sin^2(\theta_i)\right)\nonumber \\ 
&=& a_i^2 + b_i^2 \to 0 \Leftrightarrow \lim_{i \to \infty} \frac{v_i}{|v_i|} \in ker(L)\nonumber.
\end{eqnarray}

On the other hand
\begin{eqnarray}
\psi^{-1}_iv_i &=&  (\sigma_i^{-1} \circ R(\theta_i)^{-1})(v_i^\prime) \nonumber\\
&=& (\sigma_i^{-1} \circ R(-\theta_i))(a_i x_i^\prime + b_i y_i^\prime) \nonumber\\
&=& \sigma^{-1}_i\left(a_i\cos(\theta_i)x^\prime_i - a_i\sin(\theta_i)y^\prime_i\right) + \sigma^{-1}_i\left(b_i\cos(\theta_i)y^\prime_i + b_i\sin(\theta_i)x^\prime_i\right)\nonumber\\
&=& \left(a_i\cos(\theta_i) + b_i\sin(\theta_i)\right)x_i + \left(b_i \cos(\theta_i) - a_i\sin(\theta_i) \right)y_i\nonumber
\end{eqnarray}
so that
\begin{eqnarray}
|\psi_i^{-1}v_i|^2 = a_i^2 + b_i^2 \to 0 \Leftrightarrow \lim_{i \to \infty} \frac{v_i}{|v_i|}  \in halo(L)\nonumber
\end{eqnarray}
and it's evident that $|\psi_iv_i|^2 = |\psi_i^{-1}v_i|^2$. Since the last conclusions above follow from lemma~\ref{lem:asymp} we see for any sequence $v_i \in E_i$ where $v = \lim_{i \to \infty} \frac{v_i}{|v_i|}$ that $v \in ker(L) \Leftrightarrow v \in halo(L)$ when $E_i$ is a two dimensional elliptic eigenspace which persists for arbitrarily large $i$ in a sequence of symplectic maps $\left\{A_i\right\}_{i=1}^\infty$. $\qed$

We see now any sequence $\left\{A_i\right\}_{i=1}^\infty \subset Sp(2n)$ for which $Gr(A_i) \to L \in \Lambda_{2n}^n\setminus H$ has only real pairs or the usual symplectic quadruples away from the unit circle for sufficiently large $i$.$\qed$

\subsection{Decomposing Certain Unbounded Sequences of Symplectic Maps}\label{subsec:convergence}
\begin{Theorem}\label{subthm:rhothm:3}
	Any sequence $\left\{A_i\right\}_{i=1}^\infty  \subset Sp(V)$ each with unique eigenvalues for which \begin{eqnarray}
	Gr(A_i) \underset{i\to\infty}{\to} L = ker(L) \times \left\{0\right\} \oplus  \left\{0\right\} \times halo(L) \oplus Gr(\phi) \in \Lambda_{2n}^k\setminus H \nonumber
	\end{eqnarray} for some $1 \leq k \leq n-1$ with $\phi \in Sp(V_g) \cong Sp(2n-2k)$ having only semisimple eigenvalues induces a unique sequence of symplectic splittings of $V$,
	\begin{eqnarray}
	V =  E^i_s \overset{\omega}{\oplus} E^i_g\nonumber
	\end{eqnarray} into a pair of $A_i$ invariant symplectic subspaces such that $A_i = \alpha_i\oplus \beta_i \in Sp(E^i_s) \times Sp(E^i_g)$ for each $i$. Then there exists $N \in \mathbb{N}$ for which the following are true.
	\begin{enumerate}
		\item lemma~\ref{sublem:1}: There exists a subsequence for which $E^i_s \underset{i \to \infty}{\to} E_s = ker(L) \oplus halo(L) = V_s$ where $dim(E^i_s) = dim(V_s) = 2k$ and $dim(E^i_g) = 2n-2k$ for $i \in \mathbb{N}$.
		\item lemma~\ref{sublem:2}:  That $Gr(\alpha_i) \to ker(L) \times \left\{0\right\} \oplus \left\{0\right\} \times halo(L)$ and $Gr(\beta_i) \to Gr(\beta)$ for some $\beta \in Sp(E_g)$ as symplectic subspaces of $V\times \overbar{V}$.
		\item lemma~\ref{sublem:3}:  There exists a subsequence for which $E^i_g \to E_g$ and both $ker(L) \cap E^i_g = \left\{0\right\}$ and $Proj_{halo(L)}(E^i_g) = \left\{0\right\}$ for all $i \ge N$ as well as after passing to the limit. In particular observe that the second claim implies $E^i_g,E_g \leq ker(L) \oplus V_g$ for all $i \ge N$.
		\item lemma~\ref{sublem:4}: For all $i \ge N$ there exists a unique sequence of symplectic isomorphisms $I_i\colon E^i_g \to V_g$ such that $I_i = \left.Proj_{V_g}\right|_{E^i_g}$ and $I_i \underset{i\to\infty}{\to} I\colon E_g \to V_g$ where $Proj\colon dom(L) \to V_g$ is the coisotropic reduced space of $dom(L)$ identified with $V_g$.
		\item lemma~\ref{sublem:5}: Defining $\phi_i := I_i \circ \beta_i \circ I_i^{-1} \in Sp(V_g)$ then $\phi_i \to \phi\colon V_g \to V_g$ where $Gr(\phi)$ is the graph part of $L$ and the pair $\phi_i$ and $\beta_i$ share the same eigenvalues and elliptic conjugacy classes for $i \ge N$.
	\end{enumerate}
	
	From part 5 we see the pair $\phi_i$ and $\beta_i$ satisfy $\rho^2(\beta_i) = \rho^2(\phi_i)$ for $i \ge N$ leading to the main statement of the proof,
	
	\begin{eqnarray}
	\lim_{i \to \infty} \rho^2(\beta_i) = \lim_{i \to \infty} \rho^2(\phi_i) = \rho^2(\phi) =: \hat{\rho}(L)\nonumber
	\end{eqnarray} 
	
	for any sequence $\left\{A_i\right\}_{i=1}^\infty \subset Sp(2n)$ such that $Gr(A_i) \to L \in \Lambda^k_{2n} \setminus H$ for $1 \leq k \leq n - 1$.
	
\end{Theorem}

\textbf{Proof: }

\textit{Prerequisites: }We first recall that each $A \in Sp(2n)$ yields a direct sum of $V$ via symplectic generalized eigenspaces, that is 
\begin{eqnarray}
V = \bigoplus_{\lambda \in Spec(A) \cap D_2^+} E_{[\lambda]}\text{ where }D_2^+ = \left\{z \in \mathbb{C} \, | \, 0 < |z| \leq 1, Im(z) \ge 0\right\} \nonumber
\end{eqnarray} 
where $E_{[\lambda]}$ is the real eigenspace associated to the quadruple $(\lambda,\lambda^{-1},\overbar{\lambda},\overbar{\lambda}^{-1})$ satisfying $E_{[\lambda]}^\mathbb{C} = E_\lambda \oplus E_{\lambda^{-1}} \oplus E_{\overbar{\lambda}} \oplus E_{\overbar{\lambda}^{-1}} \leq V^\mathbb{C}$ and $E_\lambda$ denotes the (generalized) complex eigenspace associated to $\lambda$. Note that if $\lambda$ is an eigenvalue of $A \in Sp(V)$ then $\lambda \neq 0$, the specification that $\lambda \in D_2^+$ is simply a convenient way of picking a candidate from each quadruple as well as providing a unique limit point for unbounded eigenvalues (since the representative chosen from that quadruple tends to $0$).


 When $A$ has distinct eigenvalues this further restricts the possibilities for the above eigenspaces; We have already seen that the $E_{[\lambda]}$ for $\lambda \in (S^1 \cup \mathbb{R})\setminus \left\{0,\pm 1\right\}$ are real two dimensional symplectic subspaces on which $\left.A\right|_{E_{[\lambda]}}$ is either conjugate to a rotation by $\lambda \in S^1$ or to a hyperbolic transformation for $\lambda \in \mathbb{R}\setminus \left\{0,\pm 1\right\}$. The eigenvalue quadruples with $|\lambda| \neq 1$ and $Im(\lambda) \neq 0$ manifest as a pair of $A$ invariant real eigenspaces associated to the conjugate pairs $(\lambda,\overbar{\lambda})$ and $(\frac{1}{\lambda},\frac{1}{\overbar{\lambda}})$. A symplectic normal form for $A \in Sp(2n)$ restricted to this $4$ dimensional real vector space after picking some $\lambda = re^{i\theta}$ from a quadruple is given by the following with Darboux basis $(x,x^\prime,y,y^\prime)$,
\begin{eqnarray}
\begin{pmatrix}
R(re^{-i\theta})^{-1} & 0\nonumber \\
0 &  R(re^{-i\theta})^\textit{t}
\end{pmatrix}=
\begin{pmatrix}
\frac{1}{r}R(e^{i\theta}) & 0\nonumber \\
0 & rR(e^{i\theta})
\end{pmatrix},\, R(re^{i\theta}) = r\begin{pmatrix}
\cos(\theta) & -\sin(\theta)\\
\sin(\theta) & \cos(\theta)
\end{pmatrix}
\end{eqnarray}
so that each $A_i$ is the direct sum of a combination of the above symplectic eigenspaces.

We proceed by considering the eigenvalues of each $A_i$ for all $i \in \mathbb{N}$ as a sequence of tuples $(\lambda^i) := (\lambda^i_1,\dots,\lambda^i_{2n}) \in \mathbb{C}^{2n}$ treated as an unordered list. The space of unordered $\mathbb{C}$ tuples of length $2n$ may be identified with the orbit space $\mathbb{C}^{2n}/S_{2n}$ where $S_{2n}$ is the permutation group on $2n$ elements and the group action on $\mathbb{C}^{2n}$ is given by $\sigma((\lambda_1,\dots,\lambda_{2n})) = (\lambda_{\sigma(1)},\dots,\lambda_{\sigma(2n)})$ for any permutation $\sigma \in S_{2n}$. Following \cite{Ka} the topology induced on the space of unordered $\mathbb{C}$ tuples of length $2n$ as constructed above is identical to the one generated by the following metric,
\begin{eqnarray}
d((\lambda),(\tau)) = min_{\sigma \in S_{2n}} max_{i \leq 2n} |\lambda_{\sigma(i)} - \tau_{i}|\nonumber
\end{eqnarray}
with the helpful property that $\mathbb{C}^{2n}/S_{2n}$ is homeomorphic to $\mathbb{C}^{2n}$. 

For convenience if we impose some ordering of the $A_1$ eigenvalues we may use a recursive process to yield an \textit{essentially} unique representative for every subsequent element (since distinct permutations may both be a minimum in the above metric). Given any order for $(\lambda^1)$ we choose the order of the $i^{th}$ eigenvalue list $(\lambda^i) = (\lambda^i_1,\dots,\lambda^i_{2n})$  for any $i \ge 2$ by choosing a permutation which minimizes the above metric with respect to the previous element, i.e. $(\lambda^{i}_{\tau(j)}) = (\lambda^{i-1}_j)$ where $\tau$ is the minimizing permutation found in the above minimax metric. This allows us to treat the tuple as ordered given the order of the first. 

With this notion we let $\left\{\lambda^i_j\right\}_{j=1}^l$ denote the $l$ representatives from each eigenvalue quadruple of $A_i$ which lies in the closed upper half disc so that $V = \bigoplus_{j=1}^l E_{[\lambda^i_j]}$ for all $i\in\mathbb{N}$. As mentioned above each sequence $E_{[\lambda^i_j]}$ eventually has constant dimension for large $i$ and by compactness each possesses a limit $E_j$ for all $j\leq l$, potentially with lower dimension if distinct eigenvectors converge to each other in the limit. In our case the eigenvalues may converge but the dimension of the eigenspaces will be preserved due to the requirement that $\phi$ have semi-simple eigenvalues, precluding this possibility. \cite{Ka}
\begin{Lemma}\label{sublem:1}
	There exists a subsequence for which $E^i_s \underset{i \to \infty}{\to} E_s = ker(L) \oplus halo(L) = V_s$ where $dim(E^i_s) = dim(V_s) = 2k$ and $dim(E^i_g) = 2n-2k$ for any $i \in \mathbb{N}$.
\end{Lemma}
\textbf{Proof:} We proceed by distinguishing two possibilities for the behavior of the sequence of $A_i$ when restricted to each $E_{[\lambda^i_j]}$, either $||A||_{E_{[\lambda^i_j]}} = sup_{v \in E_{[\lambda^i_j]}} \frac{|A_iv|}{|v|} \to \infty$ or $||A||_{E_{[\lambda^i_j]}} = sup_{v \in E_{[\lambda^i_j]}} \frac{|A_iv|}{|v|} \to c_j \in \mathbb{R}$. In the first case this implies the existence of a sequence $\left\{v_i\right\}_{i=1}^\infty$ where $v_i \in E_{[\lambda^i_j]}$ for each $i$ such that $v_i \to v \neq 0$ yet $\frac{|A_iv_i|}{|v_i|} \to \infty$ so that $|A_iv_i| \to \infty$. Then by setting $v^\prime_i = \frac{v_i}{|A_iv_i|}$ then $v^\prime_i \to 0$ and we see 
\begin{eqnarray}\label{eqn:halo}
\frac{|A_iv^\prime_i|}{|v^\prime_i|} = \frac{|A_i v_i|}{|v_i|} \to \infty.
\end{eqnarray}

Thus we see that $(v^\prime_i,A_iv^\prime_i)\in Gr(\left.A_i\right|_{E_{[\lambda_j^i]}})$ for each $i$ and upon passing to the limit $(v^\prime_i,A_iv^\prime_i) = (v^\prime_i,\frac{A_iv_i}{|A_iv_i|}) \to (0,w) \in L$ where $\lim_{i\to \infty} \frac{A_iv_i}{|A_iv_i|} = w \in halo(L)$ since the $\frac{A_iv_i}{|A_iv_i|}$ belong to a compact set. Recall the above metric and ordering scheme which allows us for some fixed $j$ to identify a unique element $\lambda^i_j$ for each $i \ge 2$ so that we may form a single sequence for each of the $l$ eigenvalue quadruple representatives $\left\{\lambda^i_j\right\}_{i=1}^\infty \underset{i \to \infty}{\to} \lambda_j$ and the associated quadruple eigenspaces $\left\{E_{[\lambda^i_j]}\right\}_{i=1}^\infty$ (note that in the limit the eigenvalue is allowed to vanish). We define  
\begin{eqnarray}
S = \left\{j \leq l \, | \, \left.||A_i||\right|_{E_{[\lambda^i_j]}} \underset{i \to \infty}{\to} \infty\right\}\subset \left\{1,2,\dots,l\right\}\nonumber
\end{eqnarray}
which certainly satisfies $\left\{1,2,\dots,l\right\} = S \cup S^c$ so that we may define the symplectic $A_i$ invariant subspaces based on this condition,
\begin{eqnarray}
E^i_s := \bigoplus_{j  \in S} E_{[\lambda^i_j]},\qquad E_g^i := \bigoplus_{j \not\in S} E_{[\lambda^i_j]}. \nonumber
\end{eqnarray} 

We know that $V = E_s^i \overset{\omega}{\oplus}E_g^i $ for each $i$ and so $A_i = \alpha_i \oplus \beta_i \in Sp(E^i_s) \times Sp(E^i_g)$. Both subspaces must have constant dimension since the condition defining the two sets is binary and defined using the asymptotic behavior of the $A_i$. Since they reside in a compact space and there are no converging eigenvectors there exists a subsequence for which both converge to $E_{s}$ and $E_g$, symplectic subspaces of $V$.  For $v \in E_s$ it's true that $\lim_{i \to \infty} |A_iv| \to 0,\infty$ and for $v \in E_g$ that $lim_{i\to\infty} |A_iv| \to d \in (0,\infty)$ so that certainly $E_s \cap E_g = \left\{0\right\}$ so that $V = E_s \oplus E_g$.

\begin{Sublemma}
Equations~\eqref{eqn:squeeze} and ~\eqref{eqn:halo} imply $ker(L) \oplus halo(L) \leq E_{s}$. 
\end{Sublemma} 
\textbf{Proof:} We may construct a sequence of isotropic subspaces $\langle w^i_j \rangle_{j=1}^k \leq E^i_s$ assuming each of the $w^i_j$ converge such that 
\begin{eqnarray}
lim_{i\to\infty} \langle w^i_j \rangle_{j=1}^k = \langle w_j \rangle_{j=1}^k=halo(L)\nonumber
\end{eqnarray} 
and since each $E^i_s$ is symplectic there exists a sequence $(\tilde{w}^i_j) \leq E^i_s$ such that $\omega(w^i_j,\tilde{w}^i_j) = 1$ for all $j \leq k$ and all $i \in \mathbb{N}$. Then we consider for each $j \leq k$ the sequence of two dimensional symplectomorphisms $\left.A_i\right|_{\langle w^i_j,\tilde{w}^i_j\rangle}$ on which $|A^{-1}_iw^i_j| \to 0$. Using the argument found in lemma~\ref{lem:halo} as well as ~\eqref{eqn:squeeze} and ~\eqref{eqn:halo} the sequence $v^i_j := \frac{A^{-1}_iw^i_j}{|A^{-1}_iw^i_j|}$ satisfies $|A_iv^i_j| \to \infty$ and since $\left.\omega\right|_{\langle w^i_j,\tilde{w}^i_j\rangle \times \langle w^i_j,\tilde{w}^i_j\rangle}$ is an area form preserved by $\left.A_i\right|_{\langle w^i_j,\tilde{w}^i_j\rangle}$ for each $i$ we must have for each $j \leq k$ some sequence $\left\{z^i_j\right\}_{i=1}^\infty$ where each $z^i_j \in E^i_s$ and $z^i_j \to z_j\neq 0$ such that $|A_iz^i_j| \to 0$ for all $j \leq k$. Thus $\langle z_j \rangle_{j=1}^k = ker(L)$ implying $ker(L) \leq E_s$ so that since $ker(L) \cap halo(L) = \emptyset$ we see that $ker(L) \oplus halo(L) \leq E_s$.$\qed$

 Noting that $dim(E_g) \ge dim(V_g) = 2n-2k$ we see from above that $dim(E_s) = 2k$ and $dim(E_g) = 2n-2k$ so the corollary finishes part one of the proof. $\qed$
 
\begin{Lemma}\label{sublem:2} 
	Given the above, $Gr(\beta_i) \underset{i\to\infty}{\to} Gr(\beta)$ for $\beta \in Sp(E_g)$ and $Gr(\alpha_i) \to ker(L) \times \left\{0\right\} \oplus \left\{0\right\} \times halo(L)$ as subspaces.
\end{Lemma}

\textbf{Proof:} As we mentioned above there is a subsequence for which $E^i_g \to E_g$ so that since $Gr(\alpha_i) \oplus Gr(\beta_i) = Gr(A_i) \to L$ we see that $Gr(\beta_i) \to K_g \leq E_g \times E_g$. If $K_g \neq Gr(\beta)$ for any $\beta \in Sp(E_g)$ then $ker(K_g) \neq \left\{0\right\}$ which violates the known dimension of $ker(L)$. Indeed since $v \in ker(K_g)$ implies $v \in ker(L)$ when $K_g \leq L$ this shows that $K_g  = Gr(\beta)$ for some $\beta \in Sp(E_g)$. 

As for $Gr(\alpha_i) = \left\{(v,A_iv) \, | \, v \in V^i_s\right\}$ we may write a convergent sequence of $2k$ tuples $(v^i_j,w^i_j)_{j=1}^k \subset E^i_s$ for all $i$ such that $\langle v^i_j \rangle_{j=1}^k \to ker(L)$ and $\langle w^i_j\rangle_{j=1}^k \to halo(L)$. Then we see that $(v^i_j,A_iv^i_j) \to (v_j,0) \in ker(L) \times \left\{0\right\}$ and $(A_i^{-1}w^i_j,w^i_j) \to (0,w_j) \in \left\{0\right\} \times halo(L)$ so that indeed $Gr(\alpha_i) \to K_s = ker(L) \times \left\{0\right\} \oplus \left\{0\right\} \times halo(L)$. $\qed$

\begin{Lemma}\label{sublem:3}  
There exists an $N \in \mathbb{N}$ for which both $ker(L) \cap E^i_g = \left\{0\right\}$ and $Proj_{halo(L)}(E^i_g) = \left\{0\right\}$ for all $i \ge N$ as well as after passing to the limit.
\end{Lemma}

\textbf{Proof:} We begin by proving the following lemma,

\begin{Sublemma}
	Given $A_i$ as before with $dom(L) \pitchfork halo(L)$ and $j \leq l$ fixed then we claim $||\left.A_i\right|_{E_{[\lambda^i_j]}}||$ is bounded as $i \to \infty$ if and only if there exists an $N \in \mathbb{N}$ for which $Proj_{halo(L)} (E_{[\lambda^i_j]}) = \left\{0\right\}$ for all $i \ge N$.
\end{Sublemma}
\textbf{Proof: }First suppose $F_i := Proj_{halo(L)}(E_{[\lambda^i_j]})$ is a sequence of subspaces such that $F_i \neq \left\{0\right\}$ for all $i \in \mathbb{N}$. Then there exists a sequence $\left\{w_i\right\} \in E_{[\lambda^i_j]}$ with $w_i \to w \neq 0$ so that since $V =  halo(L)\oplus dom(L)$ we may write $w_i = f_i+g_i \to f+g = w$ where each $f_i \in halo(L)$ and $g_i \in dom(L)$. Then since $|A_i^{-1}f_i| \to 0$ we may define 
\begin{eqnarray}
v_i := \frac{A_i^{-1}(f_i +g_i)}{|A_i^{-1}f_i|}  = A_i^{-1}\left(\frac{f_i+ g_i}{|A_i^{-1}f_i|}\right) \in E_{[\lambda^i_j]}\nonumber
\end{eqnarray} and see that $|A_iv_i| = \frac{\left|f_i+ g_i\right|}{|A_i^{-1}f_i|} $. Since $g_i \to g$ with $|g| < \infty$ and the same for $f_i$ it must be $|A_iv_i| \to \infty$ and $||\left.A_i\right|_{E_{[\lambda^i_j]}}||$ is unbounded.$\qed$

Alternatively if for some $N \in \mathbb{N}$ for which $Proj_{halo(L)}(E_{{[\lambda^i_j]}}) = \left\{0\right\}$ for all $i \ge N$ we know $E_{{[\lambda^i_j]}} \leq dom(L)$ and therefore any converging sequence $v_i \in E_{{[\lambda^i_j]}}$ may be uniquely written as $v_i = k_i + g_i \in ker(L) \oplus V_g = dom(L)$ for any $i \ge N$ where $k_i+g_i \to k+g = v$. Then since $|A_ik_i| \to 0$ and $A_ig_i \to \phi(g) \in V_g$ by definition we see that the operator norm of $A_i$ over $E_{[\lambda^i_j]}$ is bounded. $\qed$

Now since $\left.||A_i||\right|_{E^i_g}$ is bounded by construction then sufficiently large $i$ and the above lemma shows that $E^i_g \leq Ker(L) \oplus V_g$. We conclude this portion with a corresponding corollary regarding the kernel.
\begin{Sublemma}
	If $||A_i||_{E_{[\lambda^i_j]}}$ is bounded as $i \to \infty$ there exists an $N \in \mathbb{N}$ for which $E_{[\lambda^i_j]} \cap ker(L) = \left\{0\right\}$ for all $i \ge N$.
\end{Sublemma}

Assume $||A_i||_{E_{[\lambda^i_j]}} \to c \in \mathbb{R}$, then $Gr(A_i) \to L$ implies that the sequence $Gr(\left.A_i\right|_{E_{[\lambda^i_j]}}) \to K \leq E_{[\lambda_j]} \times E_{[\lambda_j]}$ such that $K \leq L$ and since the norm remains bounded the limit is a graph of a function showing that $Gr(\left.A_i\right|_{E_{[\lambda^i_j]}}) \to Gr(f) \leq E_{[\lambda_j]} \times E_{[\lambda_j]}$ for some $f \in Sp(E_{[\lambda_j]})$. Since symplectic maps are non-singular we have $ker(L) \cap E_{[\lambda_j]} = \left\{0\right\}$ so that since $ker(L) \leq E_s$ and $E^i_s \oplus E^i_g = V$ we see for some $N \in \mathbb{N}$ that $ker(L) \cap E_{[\lambda^i_j]} = \left\{0\right\}$ for all $i \ge N$ as a consequence. $\qed$

\begin{Lemma}\label{sublem:4}  
	For all $i \ge N$ there exists a unique sequence of symplectic isomorphisms $I_i\colon E^i_g \to V_g$ such that $I_i = \left.Proj_{V_g}\right|_{E^i_g}$ and $I_i \underset{i\to\infty}{\to} I\colon E_g \to V_g$.
\end{Lemma}
\textbf{Proof:} Since $E^i_g$ is a $2n-2k$ dimensional subspace of $dom(L)$ with $ker(L) \cap E^i_g = \left\{0\right\}$ for sufficiently large $i$ then there eventually exists a unique symplectic map $dom(L)/ker(L) \cong E^i_g$ (refer to theorem~\ref{subthm:decomp}) for each $i$. We denote the above isomorphisms (The co-isotropic reduction of $dom(L)$ restricted to $E^i_g$) by $I_i\colon E^i_g \to dom(L)/ker(L)$ and see that $V_g \cong dom(L)/ker(L) \cong E^i_g$ uniquely for large $i$. The continuity of the coisotropic reduction with respect to a varying subspace of constant dimension shows that $I_i \to I\colon E_g \to V_g$. $\qed$

\begin{Lemma}\label{sublem:5} 
	 Defining $\phi_i := I_i \circ \beta_i \circ I_i^{-1} \in Sp(V_g)$ then $\phi_i \to \phi\colon V_g \to V_g$ where $Gr(\phi)$ is the graph part of $L$ and the pair $\phi_i$ and $\beta_i$ share the same eigenvalues and elliptic conjugacy classes for $i \ge N$.
\end{Lemma}
\textbf{Proof:} By using the above identification from section~\ref{subthm:decomp} between $V_g$ and $dom(L)/ker(L)$ we may define $\phi_i := I_i \circ \beta_i \circ I_i^{-1} \in Sp(V_g)$ so that the $\phi_i$ and $\beta_i$ are conjugate. Then since $Gr(A_i) = Gr(\alpha_i) \oplus Gr(\beta_i) \to L$ and $\beta_i \to \beta$ as $i \to \infty$ we see from the continuity of the projection and inclusion maps that 
\begin{eqnarray}
\lim_{i \to \infty} I_i \circ \beta_i \circ I_i^{-1} = I \circ \beta \circ I^{-1} = \phi \in Sp(V)\nonumber
\end{eqnarray} 
and indeed $\phi$ and $\beta$ are conjugate by $I$ as well as $\phi_i$ and $\beta_i$ via $I_i$ for sufficiently large $i$.

A subtle yet critical note here is that this notion of conjugacy occurs between distinct domains so we verify manually that $\rho^2(\phi_i) = \rho^2(\beta_i)$ for $i \ge N$. First if $\lambda_i \in Spec(\beta_i)$ with eigenvector $v^i_\lambda$ we let $w^i_\lambda = I_i(v^i_\lambda) \in V_g$ so that $(I_i \circ \beta_i \circ I_i^{-1})(w_\lambda) = (I_i \circ \beta_i)(v^i_\lambda) = I_i(\lambda_i v^i_\lambda) = \lambda_i w^i_\lambda$ and the two indeed share the same eigenvalues with $I_i$ mapping eigenvectors of $\beta_i$ to eigenvectors of $\phi_i$.

The remaining concern is regarding the preservation of the conjugacy class of the elliptic eigenvalues since they are the only eigenvalues which effect $\rho^2$. We must verify that the symmetric bilinear form written below maintains the same number of positive eigenvalues under each $I_i$, that is if $E_\lambda \leq V^\mathbb{C}$ is an elliptic eigenspace for $\beta_i$ we let
\begin{eqnarray}
Q_i\colon E_{\lambda_i} \times E_{\lambda_i} &\to& \mathbb{R}\nonumber\\
(z,z^\prime) &\mapsto& Im(\omega(z,\overbar{z^\prime})).\nonumber
\end{eqnarray} 

so that the corresponding bilinear form for $\phi_i$ will be given by $Q_i\circ (I_i^{-1} \times I_i^{-1})$ defined over the eigenspace $F_{\lambda_i} = I_i(E_{\lambda_i})$. 
\begin{eqnarray}
Q_i\circ (I_i^{-1} \times I_i^{-1})\colon F_{\lambda_i} \times F_{\lambda_i} &\to& \mathbb{R}\nonumber\\
(z,z^\prime) &\mapsto& Im(\omega(I_i^{-1}(z),\overbar{I_i^{-1}(z^\prime)})).\nonumber
\end{eqnarray}
We have implicitly extended $I_i^{-1}$ to a complex symplectic map in the natural way ($I_i^{-1}(iv) := iI_i^{-1}(v)$) so that $I_i^{-1}(\overbar{z}) = \overbar{I_i^{-1}(z)}$ meaning 
\begin{eqnarray}
Im(\omega(I_i^{-1}(z),\overbar{I_i^{-1}(z^\prime)})) = Im(\omega(I_i^{-1}(z),I_i^{-1}(\overbar{z^\prime}))).\nonumber
\end{eqnarray} 

Then since each $I_i$ is simply the co-isotropic reduction of $dom(L)$ restricted to $E_i$ we see that $\omega(I^{-1}_iz,I^{-1}_iz^\prime) = \omega(z,z^\prime)$ for any $z,z^\prime \in F^i_\lambda$ and the two bilinear forms coincide via $I_i$ and therefore share the same number of positive eigenvalues and $\rho^2(\beta_i) = \rho^2(\phi_i)$ for any $i \ge N$.
$\qed$

\section{Concluding Proofs for the Main Theorems}\label{subsec:finalproof}
\noindent\textbf{Theorem~\ref{thm:rhothm}.}
\emph{Define the map $\hat{\rho}\colon  \Lambda_{2n}\setminus H \to S^1$ by $\hat{\rho}(L) := \rho^2(\phi)$ for any $L \in \Lambda_{2n}\setminus H$ possessing the unique decomposition,}
\begin{eqnarray}
L = ker(L) \times \left\{0\right\} \oplus \left\{0\right\}\times halo(L) \oplus Gr(\phi) \leq (V_s \times V_s) \oplus (V_g \times V_g).\nonumber
\end{eqnarray} 
\emph{Then the map $\hat{\rho}$ is continuous and the diagram below commutes.}
\begin{figure}[h]
	\centering
	\begin{tikzcd}
	\Lambda_{2n} \arrow[rd,"\hat{\rho}"] \setminus H & \\
	Sp(2n) \arrow[u,hook,"Gr"] \arrow[r,"\rho^2\,\,"]& S^1
	\end{tikzcd}
\end{figure}

\textbf{Proof:} We first refer above to our implicit use of theorem~\ref{subthm:decomp} decomposing $L$ since $dom(L) \cap halo(L) = \left\{0\right\}$. For every $L = Gr(A) \in \Lambda_{2n}^0 \cong Sp(2n)$ we know $\phi = A\colon  V \circlearrowleft$ so that $\hat{\rho}(L) = \rho^2(A)$ and the above diagram commutes so that $\hat{\rho}$ extends $\rho^2$ via the graph map, it remains to show continuity.

Note that for any $L \in\Lambda^n_{2n} \cong \Lambda_n \times \Lambda_n$ that $\rho(L) = \pm 1$ for sufficiently large $i$ (\ref{subthm:hyperbolic}) so that $\hat{\rho}(L) = \lim_{i \to \infty} \rho^2(A_i) = 1$ for any sequence $A_i$ such that $Gr(A_i) \to L \in \Lambda^n_{2n}$ and $\rho^2$ may be (rather trivially) continuously extended to $\Lambda^n_{2n}$. 

For $L \in \Lambda^k_{2n}$ with $1 \leq k \leq n-1$ from theorem~\ref{subsec:convergence} we see there exists some $N \in \mathbb{N}$ for which $A_i = \alpha_i \oplus \beta_i \in Sp(E_s^i)\times Sp(E_g^i) \cong Sp(V_s) \times Sp(V_g)$ (the first coordinate of this isomorphism is arbitrary but since $\alpha_i$ diverges it is of no concern, the second coordinate isomorphism is unique for large $i$ via theorem~\ref{subsec:convergence}) such that 
\begin{eqnarray}
Gr(A_i) = Gr(\alpha_i) \oplus Gr(\beta_i) \to ker(L) \times \left\{0\right\} \oplus \left\{0\right\}\times halo(L) \oplus Gr(\beta).\nonumber
\end{eqnarray} 

For any $(v,w) \in Gr(\beta)$ we may decompose $v = v_k + v_g \in ker(L) \oplus V_g$ so that since $(v_k,0) \in ker(L) \times \left\{0\right\}$ we have that $(v_g,w) \in L$. Then since we know $Proj_{halo(L)}(E_g^i) = \left\{0\right\}$ for $i \ge N$ then $w = w_h + w_g \in halo(L) \oplus V_g$ with $w_h = 0$ so that $(v_g,w_g) \in Gr(\beta) \leq L$ since $v_g \in V_g$ means we have $w_g = \phi(v_g)$, i.e. the two graphs are seen to coincide after removing the $ker(L)$ components from the source in $Gr(\beta)$. This is simply an excessive confirmation that the normal form given in theorem~\ref{subthm:decomp} is identical to the limit of $Gr(\alpha_i) \oplus Gr(\beta_i)$ after what amounts to some column operations on $Gr(\beta)$. 

Now since $Gr(\alpha_i) \to ker(L) \times \left\{0\right\} \oplus \left\{0\right\}\times halo(L)$ with $ker(L) \cap halo(L) = \left\{0\right\}$ then theorem~\ref{subthm:hyperbolic} shows that $\rho^2(\alpha_i) = 1$ for $i\ge N$ and thus $\rho^2(A_i) = \rho^2(\alpha_i)\rho^2(\beta_i) = \rho^2(\beta_i)$ for all $i  \ge N$. Then since $\rho$ is continuous on $Sp(V_g)$ it follows that 
\begin{eqnarray}
\lim_{i \to \infty} \rho^2(A_i) = \lim_{i \to \infty}\rho^2(\beta_i) = \lim_{i \to \infty}\rho^2(\phi_i) = \rho^2(\phi) \nonumber
\end{eqnarray}
and $\hat{\rho}(L) := \rho^2(\phi)$ is indeed continuous. $\qed$

\begin{Theorem}{(theorem~\ref{thm:mainthm}, part 2)}\label{subthm:mainthm:2}
	There exists a unique real valued continuous function $\hat{\Delta}$ defined on fixed endpoint homotopy classes of paths in $\Lambda_{2n}\setminus H$ such that for any path $\gamma \in Sp(2n)$ we have that $\hat{\Delta}(Gr(\gamma)) = 2\Delta(\gamma)$.
\end{Theorem}

\textbf{Proof: }We begin by restating the construction of the extended mean index for arbitrary paths $\gamma\colon I \to \Lambda_{2n}\setminus H$ using the unique continuous map $\hat{\theta}\colon I \to \mathbb{R}$ satisfying $\hat{\theta}(0) \in [-\pi,\pi)$ and $(\hat{\rho}\circ \gamma)(t) = e^{i \hat{\theta}(t)}$ for all $t \in I$. Then if we let $\hat{\Delta}(\gamma) = \frac{\hat{\theta}(1) - \hat{\theta}(0)}{2\pi}$ it is invariant on fixed endpoint homotopy classes by construction and since $\hat{\rho} = \rho^2$ on $\Lambda_{2n}^0 \cong Sp(2n)$ its clear that for paths $\gamma \subset Sp(2n)$ that $\hat{\theta}(t) = 2\theta(t)$ so that
\begin{eqnarray}
\hat{\Delta}(Gr(\gamma)) = \frac{2\theta(1) - 2\theta(0)}{2\pi} = 2\Delta(\gamma).\nonumber\qed
\end{eqnarray}

\begin{Theorem}{(theorem~\ref{thm:mainthm}, part 3)}\label{subthm:mainthm:3}
Given $\gamma\colon I \to \Lambda_{2n} \setminus H$ and defining the product $(\gamma \cdot \gamma)(t) := \gamma(t) \circ \gamma(t)$ as canonical relations, then $\hat{\Delta}(\gamma^l) = l\cdot\hat{\Delta}(\gamma)$ for $l \ge 0$.
\end{Theorem}

\textbf{Proof:} Recalling lemma~\ref{lem:rhohomo} we see that $\hat{\rho}(\gamma^l(t)) = \hat{\rho}(\gamma(t))^l$ for all $t \in [0,1]$ and $l \in \mathbb{N}$. Then the unique angle function $\hat{\theta}$ satisfying $(\rho \circ \gamma)(t) = e^{i\hat{\theta}(t)}$ and the computation $\rho(\gamma^l(t)) = \rho(\gamma(t))^l = (e^{i\hat{\theta}(t)})^l = e^{il\hat{\theta}(t)}$ shows that $\hat{\Delta}(\gamma^l) = \frac{l\hat{\theta}(1) - l\hat{\theta}(0)}{2\pi} = l\hat{\Delta}(\gamma)$. A quick verification shows that $\gamma^0(t) = \triangle_{V}$ for any $\gamma$ so that $\Delta(\gamma^0(t)) = 0$. $\qed$

\section{Remarks}\label{sec:remarks}
\subsection{The extended Mean index on $\Lambda_2\setminus H$}\label{subsec:mean}

As mentioned above $\Delta$ may be defined on paths in $Sp(2n)$ originating at the identity so that since $\Delta$ is constant on fixed end-point homotopy classes, $\widetilde{Sp}(2n)$ may be used as the domain instead. This is defined via the association of any $\gamma$ where $\gamma(0) = Id$ with some $g \in \widetilde{Sp}(2n)$ by setting $g:= \tilde{\gamma}(1)$ (where $\tilde{\gamma}$ is the lifted path). This map is clearly onto as $\widetilde{Sp}(2n)$ is path connected and is one-to-one on fixed end-point homotopy classes of paths.

\begin{Remark}
	It is well known that $\pi_1(\Lambda_{2n}) \cong \mathbb{Z}$ but $\pi_1(\Lambda_{2n}\setminus H)$ is likely to be far larger. The same approach taken in the symplectic case to consider $\hat{\Delta}\colon  \widetilde{\Lambda_{2n}\setminus H} \to \mathbb{R}$ will no longer yield a $1-1$ map (as $\hat{\rho}$ ceases to be a circle map) meaning an intermediate covering space will be necessary to at least recover the identification (if not the group structure) which leads us to our next example.
\end{Remark}

\begin{Example}\label{example:fundamental}
	For $Sp(2) \hookrightarrow \Lambda_2$ we have a nice geometric interpretation: $Sp(2) \cong D^2 \times S^1$ and 
	\begin{eqnarray}
	\Lambda_2 = S^2\times S^1/(x,t) \sim (-x,-t) \cong \overbar{D^2}\times S^1/ \sim\nonumber
	\end{eqnarray}
	where $\sim$ identifies boundary points which are antipodal with respect to only the $S^1$ term (i.e. $(1,\theta,t) \sim (1,\theta,t+\pi)$). 
	Then 
	\begin{eqnarray}
	\Lambda_2\setminus H \cong (S^2\times S^1/\sim )\setminus \left\{(0,\theta,\pi/2)\sim (0,\theta+\pi,3\pi/2)\right\} \nonumber
	\end{eqnarray} has $\pi_1(\Lambda_2\setminus H) = \mathbb{Z}[\eta] \ast \mathbb{Z}[\tau]$ where $[\eta]$ is the push-forward of the generator for $\pi_1(Sp(2))$ and $\tau$ corresponds to a loop about the missing circle. 
	
	The minimal covering space $E \underset{Pr}{\to} \Lambda_2\setminus H$ is the one satisfying the property that $
	Pr_\ast(\pi_1(E)) = \left\{ \eta^{l}\tau\eta^{-l} \, | \, l \in \mathbb{Z} \right\} $ which is given by
	\begin{eqnarray}
	E = (\overbar{D^2}\times \mathbb{R} \underset{\phi\colon \partial \circlearrowleft}{\sqcup} \overbar{D^2}\times \mathbb{R})\setminus K \nonumber
	\end{eqnarray}
	where the set $K =  \partial D^2\times (\frac{\pi}{2}+\pi\mathbb{Z})$ and $\phi$ identifies the boundaries of the two solid cylinders via $(x,t) \sim (x,t+\pi)$. 
	\begin{Remark}
		The above subgroup of $\pi_1(\Lambda_2\setminus H)$ is in the kernel for $\hat{\Delta}$ when restricted to loops being that $\hat{\Delta}([\tau]) = 1 - 1 = 0$: it is on such a covering space for which we will regain the identification between the paths originating at the identity in $\Lambda_2\setminus H$ and points in $E$. One issue still remaining is whether $\Lambda_2 \setminus H$ possesses any richer structure (like the group structure of $\widetilde{Sp}(2n)$) so that $\hat{\Delta}$ may satisfy some analog of the algebraic properties $\Delta$ possesses.
	\end{Remark} 

\begin{figure}[h]\label{fig:graph}
	\centering
	\includegraphics[scale=0.111]{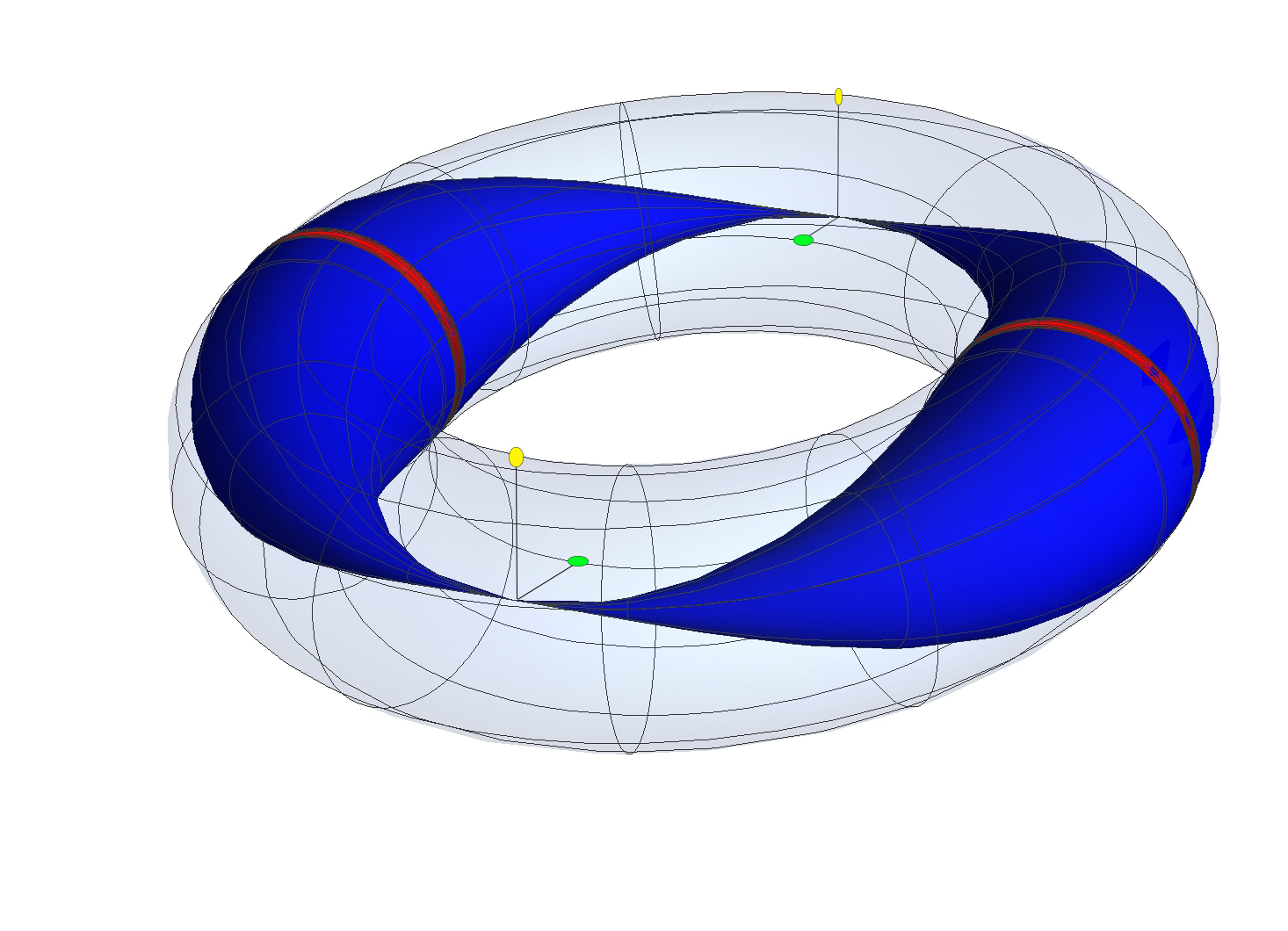}
	\caption{\emph{(Mathematica 10)} See $\Lambda_{2}^0 \cong Sp(2)$ with $H$ in red (light gray) and the two connected components of the co-oriented surface of parabolic transformations in blue (dark gray). A pair of boundary identifications are labeled by like colored ellipses. }
	\label{fig:Graph01}
\end{figure} 	
	
\end{Example}
\subsection{Smooth Canonical Relations}\label{subsec:smooth}

The extended mean index is tied significantly to linear relations, and in particular to apply this mean index to linearizations of paths along Lagrangian submanifolds will impose significant constraints on said submanifold, two of which follow below are necessary to define the mean index of a path on a Lagrangian submanifold.
\begin{enumerate}
	\item First we require $(M,\omega)$ be a symplectic manifold with $dim(M) = 4l$ for some $l \in \mathbb{N}$ and let $L$ be a Lagrangian submanifold. Then there must be a smooth bundle decomposition $TM = B \oplus C$ where $B \to M$ and $C \to M$ are $2l$ dimensional real vector bundles modeling the source and target equipped with some unique bundle isomorphism $B \cong C$  without which critical concepts such as $ker(T_xL)$ would be meaningless. It is common to consider smooth relations $M \times M$ for symplectic manifolds $M$ for which this property holds by construction.
	\item Next we would need that $T_xL \notin LagGr(T_xM) \setminus H_x$ for all $x \in L$ (or at least for all $x \in \gamma(I)$, that is the path we are linearizing must stay away from such points) where $H_x = \left\{L \in LagGr(T_xM)\, | \, \pi_{B_x}(L) \cap ( L \cap(\left\{0\right\}_x \times C_x))\right\} = \left\{0\right\}_x$ as well as some further regularity conditions as detailed in \cite{Wei} among others. This leads to the question of whether there exists any global obstructions to a Lagrangian submanifold satisfying these conditions, as well as the possibility that such Lagrangian submanifolds may be rare or non-existent for large classes of manifolds satisfying the first property. 
	\item Many of the proofs referenced in the introduction rely on $\Delta$ being a quasimorphism for which there is no obviously useful analogue in the extended case known to the author.
\end{enumerate}

The question of which smooth canonical relation framework to work in is also an immediate question. Perhaps the most promising is the Wehrheim-Woodward method applied to Lagrangian relations found in \cite{WehW} and \cite{LiW}. The \textit{highly selective category} \textbf{WW(SLREL)} consists of pairs $(L,k)$ where $L$ is a canonical relation and $k$ a non-negative integer measuring failure of transversality. The details of the construction are too nuanced for this paper but a critical step involves constructing $LagGr_\bullet(V) := \sqcup_{k = 0}^\infty LagGr(V) \times \left\{k\right\}$ equipped with a topology induced by a discrete metric which yields a weaker topology on each $LagGr(V) \times \left\{k\right\}$ than the usual one over which both the composition and reduction operations become continuous.

 Additionally the properties established regarding \textit{composable tuples} of Lagrangian relations, $(L_1,L_2,\dots,L_n)$ may in particular yield information in the context of time-dependent flows where the iterated return maps (relations) may be distinct. Similarly another benefit of the potential use of the mean index in \textbf{WW(SLREL)} is the ability to coherently form a composition of distinct $L,L^\prime\in \Lambda_{2n} \setminus H$ so that bounds of the type $|\hat{\Delta}(L \circ L^\prime) - \hat{\Delta}(L) - \hat{\Delta}(L^\prime)|$ may at least be defined, if not bounded (as is the case with $\Delta$ since it is a quasimorphism). 

Alternatively one might use the extended mean index to define the mean index of unbounded paths of symplectomorphisms converging in graph to $L \in \Lambda_{2n}\setminus H$, perhaps near unbounded punctures of pseudoholomorphic curves.

\subsection{A Proof That the Converse of Theorem~\ref{thm:rhothm}, part 1 is False}\label{thm:dimbound}
\begin{Lemma} Such a $V_g$ as described in lemma~\ref{lem:graph} with an associated unique $\phi \in Sp(V_g)$ exist if $dim((L_1^\omega + L_2^\omega)\cap (L_1\cap L_2)) \leq r  = dim(L_1^\omega \cap L_2^\omega)$ where $ Dim(L_1 \cap L_2) = 2n-2k + r$ for some $0 \leq r \leq k$. In particular this shows that the hypothesis of theorem~\ref{subthm:decomp} is not a necessary one for some $L \in \Lambda_{2n}$ to possess a uniquely determined graph portion. See \eqref{eqn:invariants} for the isotropic pair invariants used above.\end{Lemma}

\textbf{Proof:} 
We first claim that such a $V_g$ exists when
\begin{eqnarray}
dim((L_1\cap L_2)/(L_1^\omega + L_2^\omega)\cap (L_1\cap L_2)) \ge 2n-2k.\nonumber
\end{eqnarray}

It is important to note that $L_1 \cap L_2$ is no longer necessarily symplectic. The above certainly implies the existence of a $V_g \leq L_1 \cap L_2$ such that $dim(V_g) = 2n-2k$ as well as condition $(3)$ of the proof for theorem~\ref{subthm:decomp}, that is $V_g \cap L_i^\omega = \left\{0\right\}$ since $V_g \subset L_1\cap L_2$. Whether $V_g$ may be chosen to be a symplectic subspace remains to be shown. Since $(L_1\cap L_2)/((L_1^\omega + L_2^\omega)\cap (L_1\cap L_2))$ carries a unique reduced symplectic form $\omega_{red}([v],[w]) = \omega(v,w)$ we may choose symplectic $\hat{V}_g \leq (L_1\cap L_2)/((L_1^\omega + L_2^\omega)\cap (L_1\cap L_2))$ with dimension $2n-2k$ and choose $V_g \leq L_1 \cap L_2$ such that $V_g \underset{\pi_{1,2}}{\cong} \hat{V}_g$. Then since $\pi_{1,2}^\ast \omega_{red} = \omega$ on $L_1\cap L_2\leq V$ it follows that $V_g$ will be a symplectic subspace of $V$. Thus the above bound guarantees the existence of $V_g$ such that $\phi \in Sp(V_g)$.

We observe that 
\begin{eqnarray}
dim((L_1\cap L_2)/(L_1^\omega + L_2^\omega)\cap (L_1\cap L_2)) &\ge& 2n-2k \Leftrightarrow \nonumber\\
(2n-2k+r) - dim((L_1^\omega + L_2^\omega)\cap (L_1\cap L_2)) &\ge& 2n-2k \Leftrightarrow\nonumber\\
dim((L_1^\omega + L_2^\omega)\cap (L_1\cap L_2)) &\leq& r = dim(L_1^\omega \cap L_2^\omega)\nonumber.
\end{eqnarray}

In the context of theorem~\ref{thm:mainthm} we see that $\kappa = 0$ implies $r = 0$ and indeed $ dim((L_1^\omega \oplus L_2^\omega)\cap (L_1\cap L_2)) = dim((ker(L) \oplus halo(L))\cap V_g)= 0$. $\qed$

\begin{Example}\label{example:converse}
	An explicit example of the existence of Lagrangian $L$ where $\kappa(L) \neq 0$, yet $\phi$ is uniquely determined follows below. Due to the nature of the bounds imposed on the invariants found in equation~\eqref{eqn:invariants} such an $L$ exists only when $n \ge 3$, indeed since any $L \in \Lambda^1_{4}$ has either $r = 0$ or $r = 1$ then either $L \notin H$ or $dom(L) = ran(L)$ respectively. Let
	\begin{eqnarray}
	L = \langle (e_1,0),(e_2,0),(e_3,f_3),(f_3,e_3),(0,e_1),(0,f_2)\rangle \leq \mathbb{R}^6 \times \overbar{\mathbb{R}^6}\nonumber
	\end{eqnarray}
	where $(e_i,f_i)_{i=1}^3$ is a Darboux basis and $L$ is Lagrangian. 
	
	We observe that $L_1 = \langle e_1,e_2,e_3,f_3\rangle$, $L_2 = \langle e_1,e_3,f_2,f_3\rangle$, $L_1^\omega = \langle e_1,e_2\rangle$ and $L_2^\omega = \langle e_1,f_2\rangle$. Thus $dim(L_1 \cap L_2^\omega) = 1$ yet $V_g = \langle e_3,f_3\rangle \leq L_1 \cap L_2$ satisfies $V_g \cap L_i^\omega = \left\{0\right\}$ for $i =1,2$. Additionally we see that $\phi = \begin{pmatrix}
	0 & 1 \\
	1 & 0
	\end{pmatrix}$ under basis $(e_3,f_3)$. Verifying with the above lemma we see that indeed $dim((L_1^\omega + L_2^\omega) \cap (L_1 \cap L_2)) = 1 \leq r$ since $r=dim(L_1^\omega\cap L_2^\omega) = 1$.
\end{Example}

What follows is an explicit example referenced in remark~\ref{rmk:squared} regarding the failure of continuity for $\rho$ when extending to $\Lambda_{2n} \setminus H$.

\begin{Example}\label{example:squared}
	Given a Darboux basis $(e_i,f_i)_{i=1}^{n}$ for a symplectic vector space $V$ consider the two following sequences of symplectic matrices, 
	\begin{eqnarray}
	A_k &:=& Diag(\underset{\text{$n$ times}}{\underbrace{1/k,\dots,1/k}},k,\dots,k)\nonumber \\
	B_k &:=& Diag(-1/k,\dots,1/k,-k,\dots,k).\nonumber
	\end{eqnarray} 
	
	We observe for each $k$ that $\rho(A_k) = \prod_{i=1}^n 1= 1$ and $\rho(B_k) = -\prod_{i=2}^n 1 = -1$, i.e. each $A_k$ is positive hyperbolic and each $B_k$ negative hyperbolic. Then it's easy to observe that
	\begin{eqnarray}
	Gr(A_k) &=& \langle (e_1,\frac{e_1}{k}),\dots,(e_n,\frac{e_n}{k}), (f_1,kf_1),\dots, (f_n,kf_n)\rangle \nonumber\\
	&\to& \langle (e_1,0),\dots, (e_n,0), (0,f_1),\dots,(0,f_n)\rangle \nonumber\\
	Gr(B_k) &=& \langle (e_1,\frac{-e_1}{k}),\dots,(e_n,\frac{e_n}{k}), (f_1,-kf_1),\dots, (f_n,kf_n)\rangle\nonumber\\
	&\to& \langle (e_1,0),\dots, (e_n,0), (0,-f_1),\dots,(0,f_n)\rangle. \nonumber
	\end{eqnarray}
	
	Thus both $Gr(A_k) \to L\leftarrow Gr(B_k)$ yet $\rho(A_k) = 1 \neq  -1= \rho(B_k)$ while $\rho^2(A_k) = 1 = \rho^2(B_k)$ for all $k \in \mathbb{N}$.
\end{Example}

\hfill

\end{document}